\documentclass[10pt]{smfart-moi}
\usepackage{amsmath,a4wide}
\usepackage{amssymb}
\usepackage{mathrsfs}
\usepackage{euscript}
\usepackage{xspace}
\usepackage[all]{xy}
\usepackage{hyperref}
\usepackage[mathletters]{ucs}
\usepackage[T1]{fontenc}
\usepackage[utf8x]{inputenc}

\usepackage{smfthm-moi}

\usepackage[francais]{babel}
\PrerenderUnicode{éèêôÔŌōàâïî} 

\newtheorem{prop3}[paragraph]{Proposition}
\newtheorem{prp2}[subsubsection]{Proposition}
\newtheorem{thm}[subsection]{Th\'eor\`eme}
\newtheorem{thm2}[subsubsection]{Th\'eor\`eme}

\newtheorem{lmm}[subsection]{Lemme}
\newtheorem{lmm2}[subsubsection]{Lemme}
\newtheorem{lmm3}[paragraph]{Lemme}

\newtheorem{crl2}[subsubsection]{Corollaire}

\theoremstyle{definition}
\newtheorem{dfn}[subsection]{D\'efinition}

\newtheorem{rmr2}[subsubsection]{Remarque}

\newtheorem{rmrs2}[subsubsection]{Remarques}

\renewcommand{\theenumi}{\textbf{\roman{enumi}}}

\makeatletter
\renewcommand\theequation{\thesection.\textbf{\alph{equation}}}
\@addtoreset{equation}{subsection}
\makeatother

\renewcommand\thesubsubsection{\textbf{\thesubsection.\arabic{subsubsection}}} 
\renewcommand\thesubsection{\textbf{\thesection.\arabic{subsection}}}

\usepackage{amsmath}
\usepackage{amssymb}
\usepackage{mathrsfs}
\usepackage[all]{xy}

\newcommand{\Ens}{\mathsf{Ens}}

\renewcommand{\Im}{\mathrm{Im}}

\newcommand{\OO}{\mathscr{O}} 

\newcommand{\SP}{\mathrm{Spec}}

\newcommand{\et}{\mathrm{\acute{e}t}}

\newcommand{\ga}{\mathrm{Gal}}

\newcommand{\red}{\mathrm{r\acute{e}d}}

\newcommand{\dlog}{\mathrm{dlog}}
\newcommand{\cd}{{\mathrm{cd}}}

\def\car#1{\mathrm{car}.\,#1}

\newcommand{\sep}{^{\mathrm{s\acute{e}p}}}

\newcommand{\RG}{\mathrm{R\Gamma}}

\newcommand{\HH}{\mathrm{H}}

\newcommand{\TR}{\mathrm{Tr}}

\renewcommand{\ker}{\mathrm{Ker}}

\newcommand{\rang}{{\mathrm{rang}}}
\newcommand{\Frac}{{\mathrm{Frac}}}
\newcommand{\prang}{{p\text{-}\mathrm{rang}}}

\newcommand{\ZZ}{\mathbf{Z}}

\newcommand{\NN}{\mathbf{N}}

\newcommand{\QQ}{\mathbf{Q}}

\newcommand{\FF}{\mathbf{F}}

\newcommand{\MM}{\mathfrak{m}}

\newcommand{\ra}{\rightarrow}

\newcommand{\hra}{\hookrightarrow }
\newcommand{\inj}{\hookrightarrow}
\newcommand{\sr}{\stackrel}

\newcommand{\rrraxy}{\ar@<1ex>[r] \ar@<-1ex>[r] \ar[r] }
\newcommand{\iso}{\stackrel{\sim}{\ra}}

\newcommand{\isononcan}{\simeq}

\newcommand{\surj}{\twoheadrightarrow}

\newcommand{\cad}{c'est-\`a-dire\xspace}

\newcommand{\GG}{\mathbf{G}}

\newcommand{\colim}{{\mathrm{colim}}}

\def\sga#1#2#3{[{\bf $\mathbf{SGA_{#1}}$}~{\sc #2}~#3]}
\def\ega#1#2{[{\bf ÉGA}~{\sc #1}~#2]}
\def\ac#1#2#3#4{[{\bf Bourbaki}, A.C. chap.~{\sc #1}, §#2, N°#3\,#4]}
\def\]{\textup{\mbox{]\hspace{-.15em}]}}}
\def\[{\textup{\mbox{[\hspace{-.15em}[}}}
\def\mc{\mathscr}
\def\got{\mathfrak}

\def\-l{\ \vspace{-4mm}}

\def\L#1{\mathsf{L}_{#1}}
\def\chap#1{\widehat{#1}}

\def\sur{\overline}
\def\sous{\underline}
\def\gtilde{\widetilde}

\usepackage{textcomp}

\author{Ofer Gabber}
\address{CNRS et IHÉS\\
Le Bois-Marie\\
35, route de Chartres\\
\oldstylenums{91440} Bures-sur-Yvette\\
France}
\email{Gabber@ihes.fr}

\author{Fabrice Orgogozo}
\address{CNRS et Centre de math\'ematiques Laurent Schwartz\\
\'Ecole polytechnique\\
\oldstylenums{91128} Palaiseau\\
France}
\email{Fabrice.Orgogozo@math.polytechnique.fr}
\urladdr{http://www.math.polytechnique.fr/~orgogozo/}

\title{Sur la $p$-dimension des corps}

\begin{document}

\begin{center}
Sur la $p$-dimension des corps.\\
Ofer Gabber \& Fabrice Orgogozo\\ 2007-5-7
\end{center}

\section{Introduction}

Soient $k$ un corps et $p$ un nombre premier. Rappelons la
définition de la $p$-dimension de $k$ (\cite{CG@Serre}, \cite{Galois@Kato}).

Si $p≠\car{k}$, on appelle \emph{$p$-dimension} de $k$, 
la $p$-dimension cohomologique du 
groupe profini $G_k:=\ga(k^{\sep}/k)$, 
où $k^{\sep}$ est une clôture séparable de $k$ (cf. \cite{CG@Serre}). 
C'est le plus petit entier $d$ (où l'infini si un tel entier n'existe pas)
tel que pour tout $G_k$-module discret $M$ de $p$-torsion
et tout $n>d$, les groupes $\HH^n(G_k,M)$ soient nuls. 

Si $p=\car{k}$, la définition fait intervenir des invariants différentiels.
Pour tout schéma $X$ et tout entier $i∈\NN$, notons 
$Ω^i_X:=\bigwedge^i Ω¹_{X/\ZZ}$ le $\OO_X$-module des $i$-formes
différentielles absolues et $Ω^i_{X,\log}$ le sous-faisceau étale abélien 
des formes différentielles \emph{logarithmiques}, \cad localement engendré
par les sections de la forme $\dlog(x₁)\wedge \cdots \wedge \dlog(x_i)$, pour
$x₁,\dots,x_i∈\OO_X^×$. On pose alors 
$\HH^i_p(k):=\HH^{i}_{\et}(k,Ω^{i-1}_{k,\log}[-(i-1)])$ ;
c'est un analogue de $\HH^{i}_{\et}(k,μ_p^{\otimes i-1})$
en caractéristique différente de $p$. (Voir \ref{Hrp} pour
une définition équivalente.)
Enfin, rappelons que le rang du $k$-module $Ω¹_k$ est égal au 
\emph{$p$-rang} de $k$, \cad au cardinal
d'une $p$-base (absolue) de $k$ (cf. \ega{iv}{chap. 0 §21.1}) ; 
s'il est fini c'est également l'entier $r$ pour lequel $[k:k^p]=p^r$.

\begin{dfn}[Kazuya Katô, \cite{Galois@Kato}, §0]
La \emph{$p$-dimension} d'un corps $k$ de caractéristique $p>0$ 
est le plus petit
entier $d$ (où l'infini si un tel entier n'existe pas)
tel que $Ω^{d+1}_k=0$ et $\HH^{d+1}_p(k')=0$ pour toute extension
finie $k'/k$.
\end{dfn}

L'objet de cet article est de démontrer le résultat suivant, conjecturé
par K. Katô :

\begin{thm}\label{énoncé}
Soit $A$ un anneau local hensélien excellent, intègre de dimension $d$.
Soient $k$ son corps résiduel, de caractéristique
$p>0$, et $K$ son corps des fractions. Alors, on a l'égalité
$$
\dim_p(K)=\dim(A)+\dim_p(k).
$$
\end{thm}

Dans \emph{op. cit.}, ce théorème est démontré par K. Katô 
dans le cas particulier essentiel où $A$ 
est un anneau de valuation discrète complet (cf. \ref{kato}).
Son théorème est une généralisation d'un théorème de S. Lang 
(\cite{CG@Serre}, chap. II, §3.3 et §4.3) au cas d'un corps résiduel 
non algébriquement clos. La nécessité de prendre en compte le 
$p$-rang dans le cas d'un corps résiduel non nécessairement parfait
avait été conjecturée par M. Artin dans \sga{4}{x}{2.2}.
La démonstration de K. Katô utilise la $K$-théorie de Milnor, qui permet
en caractéristique mixte,
par l'intermédiaire des symboles cohomologiques, différentiels,
et du théorème de Bass-Tate (\cite{Bass-Tate}, I 4.3), 
de faire le pont entre la cohomologie galoisienne du corps des fractions
et les formes différentielles absolue sur le corps résiduel.

Le cas de la dimension deux est également établi par K. Katô,
en caractéristique mixte (et dans le cas d'un corps résiduel algébriquement
clos), dans \cite{Arithmetic@Saito}, §5. Sa démonstration, $K$-théorique,
repose sur le théorème de Merkurjev-Suslin ainsi
que sur la résolution des singularités des surfaces.

Notons que si A est un anneau strictement local
excellent, intègre, de corps des fractions $K$,  et
$\ell$ un nombre premier \emph{inversible dans $A$}, on a
$$\dim_\ell\,K = \dim\,A,$$
comme conjecturé par M. Artin dans \sga{4}{X}{3.1}. On a
en effet $\dim_\ell\, K ≥ \dim \, A$ (\sga{4}{X}{2.4}) ; par
ailleurs, le premier auteur 
a récemment démontré que pour tout $f ∈ A$, on a $\mathrm{cd}_\ell(\SP(A[f^{-1}]) 
≤ \dim\, A$ (\cite{Conference-Deligne@Gabber}, §8) et donc $\dim_\ell \, K ≤ \dim\,A$  
par passage à la limite.

Donnons brièvement quelques indications sur la méthode utilisée ici,
qui suit de près la technique d'algébrisation introduite par
le premier auteur dans \cite{Conference-Illusie@Gabber} et
\emph{op. cit.} (voir également \cite{Hasse@Matsumi}, théorème 2.2). 
Utilisant le théorème 
d'approximation de Popescu, on se ramène au cas d'un anneau local complet
noethérien. En égale caractéristique, on utilise alors le théorème
de Cohen-Gabber (\cite{Conference-Illusie@Gabber}, 8.1), dont la démonstration 
est rappelée en appendice (cf. \ref{gabber}), qui précise
le théorème de structure de Cohen et permet, grâce au 
théorème d'algébrisation d'Elkik, de faire de cet anneau
le complété d'un anneau local hensélien essentiellement de type fini et 
de dimension relative un sur un anneau local complet de dimension un de moins. 
Les prémices d'une
telle idée se trouvent déjà dans l'exposé de M. Artin \sga{4}{xix}{§1 \& §6}.
En caractéristique mixte, dans le cas où $p$ est ramifié dans $A$,
on utilise le théorème de Epp,
ainsi que le théorème \ref{énoncé} en égale caractéristique, pour pouvoir
algébriser nos données.  Dans les deux cas, on procède par récurrence
sur la dimension de l'anneau, 
en utilisant le théorème de Katô (dimension un) et, en caractéristique mixte, 
un théorème de comparaison hensélien/formel dû à K. Fujiwara et au premier
auteur.

Le plan de l'article est le suivant : après quelques rappels et compléments
sur la $p$-dimension (§\ref{rappels}), on commence par minorer la $p$-dimension du corps
des fractions (§\ref{section minoration}). La démonstration est très semblable à celle de K. Katô
en dimension deux. Majorer la $p$-dimension est plus difficile. On commence
par le cas d'égale caractéristique (§\ref{section majoration egale}),
qui nous permet de traiter ensuite 
le cas d'inégale caractéristique (§\ref{section majoration mixte}). 
En inégale caractéristique, le théorème principal est généralisé
(§\ref{generalisation}) au cas d'un ouvert affine de $\SP(A[p^{-1}])$.
Enfin, dans un appendice (§\ref{appendice}), on rappelle
la démonstration du théorème de Cohen-Gabber mentionné ci-dessus.

\vskip1cm
Le second auteur souhaite remercier Luc Illusie pour d'utiles
commentaires sur les versions précédentes de ce texte, ainsi
que Takeshi Saitô et l'université de Tôkyô pour leur
chaleureux accueil durant le premier semestre 2006-2007.

\section{$p$-dimension : rappels et compléments}\label{rappels}

Dans cette section on réunit divers lemmes (dont certains
ne sont mis que pour mémoire) qui seront 
utiles aux cours des dévissages qui vont suivre (réduction au cas normal, resp.
complet) ainsi que l'énoncé du théorème de Katô et d'un corollaire
important.

\subsection{$p$-rang}

\begin{lmm2}\label{invariance p-rang}
Soit $k'/k$ une extension finie de corps de caractéristique $p>0$.
Alors, le $p$-rang de $k$ est égal au $p$-rang de $k'$.
\end{lmm2}

\begin{proof}
On se ramène immédiatement aux cas où $k'/k$ est étale ou radicielle
de degré $p$.
Dans le premier cas, le morphisme canonique $k'\otimes_k Ω¹_k → Ω¹_{k'}$
est un isomorphisme et le résultat s'obtient alors en considérant le rang
de ces modules. Dans le second cas, considérons $a∈k'$ tel que $k'=k(a)$ et
posons $b:=a^p∈k$. L'élément $b$ n'appartient pas à $k^p$ et l'on peut donc 
le compléter en une $p$-base $\{b,b_i\}_{i∈I}$ de $k$. 
De l'égalité ${k'}^p=k^p(b)$, on déduit que
l'ensemble $\{a,b_i\}_{i∈I}$ forme une $p$-base de $k'$.
\end{proof}

\begin{lmm2}\label{p-rang par complétion}
Soit $K$ un corps discrètement valué de caractéristique $p$, de complété $\chap{K}$.
Si $[K:K^p]$ est fini, on a l'inégalité :
$$
[K:K^p]≥[\chap{K}:\chap{K}^p].
$$
\end{lmm2}

\begin{proof}
Soit $\{b_1,\dots,b_r\}$ une $p$-base de $K$, de sorte qu'en particulier
$K=K^p(b₁,\dots,b_r)$. Le sous-corps $\chap{K}^p(b₁,\dots,b_r)$ de $\chap{K}$
contient donc $K$ et, étant fini sur $\chap{K}^p$, est fermé dans $\chap{K}$.
L'égalité $\chap{K}=\chap{K}^p(b₁,\dots,b_r)$, et la conclusion, en découlent.
\end{proof}

\begin{rmrs2}
a) Il se peut que le $p$-rang de $K$ soit dénombrable
et celui de $\chap{K}$ indénombrable (cf. \cite{Imperfection@Bastos},~§3).
b) Le lemme est valable sans restriction sur la valuation : 
d'après \ac{6}{8}{2}{ Prop. 2}, le carré commutatif 
$$
\xymatrix{
K \ar[r] & \chap{K} \\
K \ar[u]^{\mathrm{Frob}} \ar[r] & \chap{K}
\ar[u]_{\mathrm{Frob}}
}
$$
induit une surjection $\chap{K}⊗_K K \surj \chap{K}$. L'inégalité
désirée en résulte immédiatement.
\end{rmrs2}

Réciproquement :

\begin{lmm2}\label{p-rang complété excellent}
Soit $A$ un anneau local hensélien \emph{excellent} intègre de caractéristique
$p>0$, de corps des fractions
$K$. Soient $\chap{A}$ le complété de $A$ et $\chap{K}$ son corps
des fractions. Alors, 
$$
[K:K^p]≤[\chap{K}:\chap{K}^p].
$$
De plus, si $[K:K^p]$ est fini, c'est une égalité.
\end{lmm2}

Rappelons que $\chap{A}$ est \emph{intègre}
(\ega{$\textsc{iv}_{4}$}{18.9.2}).

\begin{démo}
L'extension $\chap{K}/K$ est séparable de sorte que
le morphisme canonique $\chap{K}\otimes_K  Ω¹_K→
Ω¹_{\chap{K}}$ est une injection (\ega{$0_{\textrm{IV}}$}{20.6.3}).
Vérifions la seconde assertion. Sous les hypothèses faites, la normalisation
de $A$ dans $K^{1/p}$ est finie sur $A$, de sorte que
$\mathrm{Frob}:A→A$
est fini et que le morphisme $\chap{A}⊗_{A,\mathrm{Frob}} A→ \chap{A}$ 
est un isomorphisme (cf. p. ex. \ega{$0_{\textrm{I}}$}{7.3.3} ; rappelons à cette
occasion qu'un anneau excellent est noethérien). 
On en tire immédiatement que $Ω¹_{A}⊗_A \chap{A}→
Ω¹_{\chap{A}}$ est un isomorphisme. (Cela résulte du fait
que $Ω¹_{A}=Ω¹_{\mathrm{Frob}:A→A}$ et de même pour $\chap{A}$.)

\end{démo}


\begin{lmm2}\label{p-rang séries formelles}
Soit $A$ un anneau de caractéristique $p>0$ possédant
une $p$-base \emph{finie} $\{b_i\}_{i∈I}$. Pour tout entier $n≥0$, 
l'ensemble $\{b_i\}_{i∈I}\cup \{x_1,\dots,x_n\}$ constitue
une $p$-base de l'anneau $A[[x₁,\dots,x_n]]$.
\end{lmm2}

\begin{démo}
Il suffit de traiter le cas où $n=1$.
Pour $θ:I→[0,p-1]$, posons $b^θ=\prod_{i\in I} b_i^{θ(i)}$.
La conclusion résulte des décompositions :
$$A[[X]]=\bigoplus_{i=0}^{p-1} A[[X^p]] X^i,$$
$$
A=\bigoplus_{θ∈[0,p-1]^{I}} A^p b^θ ,
$$
et de la finitude de l'ensemble $I$.
\end{démo}

Enfin, signalons le lemme suivant : 

\begin{lmm2}\label{p-rang fractions}
Soit $A$ un anneau intègre ayant une $p$-base $\{b_i\}_{i∈I}$.
Alors, les éléments $b_i$ forment une base de $\Frac{A}$.
\end{lmm2}

\begin{démo}
C'est immédiat en chassant les dénominateurs par une puissance
$p$-ième.
\end{démo}

\subsection{Les groupes $\HH^i_p(k)$}

\begin{prp2}\label{Hrp}
Soient $k$ un corps et $i$ un entier.
Le groupe $\HH^{i+1}_p(k):=\HH^1_{\et}(k,Ω^{i}_{k,\log})$
est isomorphe au conoyau du morphisme 
$$
℘=1-C^{-1}: Ω^i_k \ra Ω^i_k/dΩ^{i-1}_{k}
$$
où $1$ est la projection canonique $Ω^i_k→ Ω^i_k/dΩ^{i-1}_{k}$
et $C^{-1}$ l'isomorphisme de Cartier inverse.
Le morphisme $℘$ est caractérisé par la formule :
$$
x \,\dlog(y₁)\wedge \cdots \dlog(y_i) \mapsto (x-x^p)\, \dlog(y₁)\wedge \cdots
\wedge \dlog(y_i) \ \mathrm{mod}\ dΩ^{i-1}_{k}.
$$
\end{prp2}

\begin{proof}
Rappelons que l'on a une suite exacte (de faisceaux étales abéliens sur $\SP(k)$)
$$
0→Ω^{i}_{k,\log}→({Ω^i_k})_{d=0}\sr{1-C}{→}Ω^i_k→0, 
$$
où $C$ est l'opérateur de Cartier sur les formes fermées (cf.
\cite{Rham-Witt@Illusie} chap. $0$, §2.4 et \cite{Syntomic@Tsuji}, 6.1.1.). 
(Pour $i=1$, l'exactitude à gauche est la caractérisation classique due à P. Cartier
des $1$-formes logarithmiques.)
On déduit immédiatement de l'isomorphisme de Cartier que le noyau de $1-C$
s'identifie naturellement au noyau de $℘$.
Il suffit alors de considérer la suite exacte longue de cohomologie associée et 
l'acyclicité des faisceaux cohérents pour conclure.
(Voir aussi \cite{Swan@Kato}, §1.3.)
\end{proof}

\subsubsection{Exemples}Pour tout corps $k$ de caractéristique $p>0$, on a
$\HH¹_p(k)=k/℘(k)$, où $℘$ est le morphisme
d'Artin-Schreier usuel. En particulier, il est trivial
pour un corps séparablement clos. Si $k$ est parfait
(de façon équivalente : de $p$-rang nul), $\HH^2_p(k)$
--- qui est en général un quotient de $Ω¹_k$ --- est nul.
On peut montrer plus précisément
que $\HH^2_p(k)$
s'identifie au sous-groupe de $\mathrm{Br}(k)$ formé
des éléments annulés par $p$, par l'intermédiaire de l'application envoyant une forme
différentielle $ω=x\dlog(y)$ ($x∈ k$, $y∈k^×$), 
sur l'algèbre centrale simple de rang $p^2$ 
définie par des générateurs $X,Y$ liés par les relations 
$X^p-X=x$, $Y^p=y$ et $YXY^{-1}=X+1$ (cf.
\cite{CG@Serre}, chap. \textsc{ii}, §2.2 et
\cite{CL@Serre}, chap. \textsc{xiv}, §5).

\subsubsection{Trace}Remarquons maintenant que l'on peut définir une trace dans le
présent contexte. La fonctorialité \emph{contravariante} en le corps est quant à
elle élémentaire. 

Si $k'/k$ est une extension finie étale de corps de caractéristique $p>0$, la
trace $\TR^{i,Ω}_{k'/k}: Ω^i_{k'}→Ω^i_k$ (déduite du morphisme $\TR_{k'/k}:k'→k$ 
grâce à l'isomorphisme $k'\otimes_k Ω^i_k \iso Ω^i_k$) 
envoie $dΩ^{i-1}_{k'}$ dans $dΩ^{i-1}_k$
et induit un morphisme, noté $\TR^{i,\HH}_{k'/k}$, de $\HH^{i+1}_p(k')$
dans $\HH^{i+1}_p(k)$. 

Si $k'/k$ est finie non nécessairement étale, de caractéristique $p>0$,
une telle trace est construite par K. Katô à partir de la norme 
$\mathrm{N}^K_{k'/k}$ en $K$-théorie de Milnor 
(cf. \cite{Bloch-Kato}, p. 126 ; voir également \cite{Explicit@Fukaya}, §2). 
Cette trace est caractérisée par les propriétés suivantes :
\begin{enumerate}
\item pour $i=0$, c'est la trace usuelle,
\item le diagramme suivant est commutatif :
$$
\xymatrix{
K_i^{\mathrm{M}}(k') \ar[d]_{\mathrm{N}^{i,K}_{k'/k}} \ar[r]^{\mathrm{dlog}} &
\Omega^i_{k'}
\ar[d]^{\mathrm{Tr}^{i,Ω}_{k'/k}} \\
K_i^{\mathrm{M}}(k) \ar[r]^{\mathrm{dlog}} & \Omega^i_k 
}
,$$
\item on a compatibilité avec $d$ :
$d(\mathrm{Tr}^{i,Ω}_{k'/k}(ω))=\mathrm{Tr}^{i+1,Ω}_{k'/k}(dω)$ (pour
$ω∈Ω_{k'}^{i}$),
\item elle satisfait une formule de projection : $\mathrm{Tr}^{i+j}(ω\wedge ω')
=\mathrm{Tr}^{i}(ω)\wedge ω'$ si $ω∈Ω^i_{k'}$ et $ω'∈Ω^j_{k}$,
\item elle est compatible aux extensions de corps : $\mathrm{Tr}^{i,Ω}_{k'/k}
\circ \mathrm{Tr}^{i,Ω}_{k''/k'}=\mathrm{Tr}^{i,Ω}_{k''/k}$.
\end{enumerate}

\label{rem trace}
Si $k'=k(a)/k$, où $a^p=b∈k$, est une extension radicielle de degré $p$, 
il résulte de la commutativité du diagramme ci-dessus
que l'on a $\mathrm{Tr}^{1,Ω}_{k'/k}(\dlog(a))=\dlog(b)$, et de
la formule de projection que l'on a, pour tout entier $i≥0$ :
$$
\TR^{i,Ω}_{k'/k}\Big((\sum_{j=0}^{p-1} c_j a^j)\, \dlog(a)\wedge \dlog(b₁) 
\wedge \cdots \wedge \dlog(b_{i-1})\Big)=
c_0 \, \dlog(b)\wedge \dlog(b_1)\wedge \cdots \wedge \dlog(b_{i-1}),
$$
où les $c_j$ appartiennent à $k$ et les $b_j$ à $k^×$. 

On déduit de cette formule que la trace commute au morphisme 
$C^{-1}:Ω^i→ Ω^i/dΩ^{i-1}$\footnote{Dans le
cas radiciel de degré $p$, on se ramène aisément au fait
que pour chaque $0<j<p$ et tous $c∈k$,
$b,b'₁,\dots,b'_{i-1}∈k^×$, 
on a $c^p b^{j-1}db\wedge
\dlog(b')=d(\frac{c^p}{j}b^j \dlog(b'))∈dΩ^{i-1}_k$,
où l'on pose $\dlog(b'):=\dlog(b'₁)\wedge\cdots \wedge \dlog(b'_{i-1})$.
Dans le cas étale, cela résulte de ce que $\mathrm{Tr}(a)^p=\mathrm{Tr}(a^p)$.} 
de sorte qu'elle induit 
par passage au quotient un morphisme trace 
$\TR^{i,\HH}_{k'/k}$ sur les $\HH^{i+1}_p$. De plus, 
on constate que le morphisme induit est une \emph{surjection} pour 
$i$ égal au $p$-rang de $k$.

Le lemme immédiat suivant rend plus explicite la structure de $\HH^{r+1}_p(k)$,
où $r$ est le $p$-rang d'un corps $k$ comme quotient de $Ω^{r}_k\isononcan k$.
(Voir aussi \ref{variante tous degrés} pour une démonstration
du lemme \ref{Hr complété} qui n'en fait pas usage.)

\begin{lmm2}\label{structure Hr}
Soit $k$ un corps de caractéristique $p>0$ et $b_1,\dots,b_r$ une $p$-base.
Pour toute fonction $θ:[1,r]→[0,\dots,p-1]$, posons 
$b^θ:=b_1^{θ(1)}\cdots b_r^{θ(r)}$ et notons $k_{>0}$ le
sous-$k^p$-espace vectoriel $\oplus_{θ≠0} k^p b^θ$ de $k$.
L'isomorphisme $k→Ω^r_k$, défini par $λ\mapsto λ \dlog(b):=
λ \dlog(b₁)\wedge \cdots \wedge \dlog(b_r)$ induit par passage au quotient
un isomorphisme de $\FF_p$-espaces vectoriels
$$
k/\big(℘(k) + k_{>0} \big) \iso \HH^{r+1}_p(k).
$$
\end{lmm2}

\begin{démo}
Calculons l'image $dΩ^{r-1}_k\subset Ω^r_k=k \cdot \dlog(b)=\oplus_{θ} 
\big( k^p \cdot b^θ \dlog(b)\big)$.
Pour chaque $θ,i$, considérons l'élément 
$$ω_{θ,i}:=b^θ \dlog(b₁)\wedge \cdots \wedge \dlog(b_{i-1}) \wedge \dlog(b_{i+1})\wedge
\cdots \wedge \dlog(b_r)$$ de $dΩ^{r-1}_k$. 
Ces formes engendrent $Ω^{r-1}_k$ comme $k^p$-espace vectoriel.
On a $$dω_{θ,i}=(-1)^{i+1} θ(i) b^θ\dlog(b)=\big((-1)^{i+1}θ(i)\big)^p b^θ\dlog(b).$$
Ainsi, faisant varier $θ,i$, on en déduit que 
$$
dΩ^{r-1}_k=k_{>0} \cdot \dlog(b).
$$

Le quotient $Ω^r_k/dΩ^{r-1}_k$ s'identifie donc
à $k/k_{>0}$ par l'application $λ \mod k_{>0}
\mapsto λ \dlog(b) \, \mod dΩ^{r-1}_k$,
et le morphisme $℘:Ω^r_k→Ω^r_k/dΩ^{r-1}_k$ à 
$λ \mapsto λ-λ^p \, \mod k_{>0}$. Ainsi, $\HH^{r+1}_p(k)$
est isomorphe à
$$
k/\big(k_{>0}+℘(k)\big).
$$
\end{démo}

Voici maintenant un analogue du lemme \ref{invariance p-rang}.

\begin{lmm2}\label{Hr complété}
Soit $K$ un corps discrètement valué de caractéristique $p>0$,
de $p$-rang $r<+∞$. L'application canonique
$$
\HH^{r+1}_p(K)→\HH^{r+1}_p(\chap{K})
$$
est une surjection.
En particulier, si $\HH^{r+1}_p(K)$ est nul, il en 
est de même de $\HH^{r+1}_p(\chap{K})$.
\end{lmm2}

\begin{démo}
Si le $p$-rang de $\chap{K}$ est strictement inférieur à $r$, 
$\HH^{r+1}_p(\chap{K})$ est nul et il n'y a rien à démontrer.
Supposons le donc égal à $r$. Le morphisme déduit de la fonctorialité
contravariante en le corps correspond, d'après le lemme \ref{structure Hr}
et sa démonstration, au morphisme canonique 
$$
K/\big( ℘(K)+K_{>0}\big) → \chap{K}/\big(℘(\chap{K})+\chap{K}_{>0}\big).
$$
(Les choix de $K_{>0}$ et $\chap{K}_{>0}$ se font relativement à une $p$-base
commune.)
Il nous faut suffit donc de montrer que le morphisme composé
$$
K\inj\chap{K}\surj\chap{K}/\big(℘(\chap{K})+\chap{K}_{>0}\big)
$$
est une surjection.
On va montrer plus précisément que le morphisme
$K→\chap{K}/℘(\chap{K})$ est une surjection.
Soit $λ∈\chap{K}$ ; considérons $λ₀∈K$ tel que $v(λ-λ₀)>0$.
D'après le lemme ci-dessous, il existe $α∈\chap{K}$
tel que $λ-λ₀=℘(α)$. Ainsi, $λ \equiv λ₀$ modulo $℘(\chap{K})$.
\end{démo}

\begin{lmm2}\label{surjectivité P}
Pour tout  anneau local $A$, complet de caractéristique $p>0$,
le morphisme $℘:A→A$, $a\mapsto a-a^p$, induit 
une surjection $\MM_A→\MM_A$.
\end{lmm2}

\begin{démo}
Soit $a∈A$ ; l'identité $a=℘(a)+a^p$ entraîne par récurrence l'égalité
$$a=℘(a+a^p+\cdots+a^{p^n})+a^{p^{n+1}}$$ pour tout $n≥0$.
Pour $a∈\MM_A$, la suite $(a+a^p+\cdots+a^{p^n})_{n≥0}$ converge 
dans $A$ vers un élément $b∈\MM_A$ ; d'après l'égalité ci-dessus on
a alors $a=℘(b)$.
\end{démo}

\begin{rmr2}\label{variante tous degrés}
La conclusion du lemme \ref{Hr complété} est vraie pour tous les
groupes $\HH_p^{i+1}$ (et pas seulement pour $i=r$). Il suffit 
pour cela de montrer que $Ω^i_{\chap{K}}$ est engendré (comme
groupe abélien) par 
$Ω^i_K$ et les formes $α \frac{\mathrm{d}b₁}{b₁}\wedge \cdots \wedge \frac{\mathrm{d}b_i}{b_i}$,
où $α∈\MM_{\OO_{\chap{K}}}$ (l'idéal maximal de l'anneau des entiers) et
les $b_j$ sont dans $\chap{K}^×$. Étant donné une forme $β 
\frac{\mathrm{d}x₁}{x₁}\wedge \cdots \wedge \frac{\mathrm{d}x_i}{x_i}$, 
$β,x_j∈\chap{K}^{×}$,
on peut trouver des $y_j∈K$ tels que $v(x_j-y_j)>v(x_j)$ et
$v(x_j-y_j)+v(β)>v(x_j)$. Posons $α_j:=\frac{x_j}{y_j}-1$,
de sorte que $x_j=y_j(α_j+1)$ et donc (pour chaque $j$)
$$
\frac{\mathrm{d}x_j}{x_j}=\frac{\mathrm{d}y_j}{y_j}+\frac{α_j}{α_j+1}
\frac{\mathrm{d}α_j}{α_j}.
$$
L'hypothèse sur les valuations signifie que $\frac{α_j}{α_j+1}$ et
$β\frac{α_j}{α_j+1}$ appartiennent à
$\MM_{\OO_{\chap{K}}}$.

Cela nous permet de remplacer les $x_j$ par les $y_j$ ; on conclut
en approchant $β$ comme ci-dessus. 
\end{rmr2}

En dimension supérieure on a la réciproque suivante : 

\begin{prp2}\label{Hr complété excellent}
Soit $A$ un anneau local hensélien excellent intègre de corps des
fractions $K$. Soient $\chap{A}$ son complété et $\chap{K}$ le corps
des fractions de $\chap{A}$. Supposons les $p$-rangs de $K$ et $\chap{K}$
égaux à un entier $r<+∞$.
Le morphisme canonique
$$
\HH^{r+1}_p(K)→\HH^{r+1}_p(\chap{K})
$$
est une injection.
\end{prp2}


La démonstration fait usage de la généralisation suivante du théorème
d'approximation d'Artin :

\begin{thm2}[Dorin Popescu, \cite{Neron@Popescu}, théorème 1.3]\label{popescu}
Soit $A$ un anneau local excellent hensélien. Pour tout système système fini
d'équations polynomiales à coefficients dans $A$, l'ensemble
des $A$-points est \emph{dense}, pour la topologie $\MM_A$-adique,
dans l'ensemble des $\chap{A}$-points.
\end{thm2}

\begin{démo}[Démonstration de la proposition \ref{Hr complété excellent}]
Soit $\{b₁,\dots,b_r\}$ une $p$-base de $K$ ; c'est également
une $p$-base de $\chap{K}$ (\ref{p-rang complété excellent}). 
Quitte à les multiplier par une puissance
$p$-ième convenable, on peut supposer les $b_i$ dans $A$. 
Soient $K_{>0}$ et $\chap{K}_{>0}$ comme
en \ref{structure Hr}, relativement à ces bases.
Il s'agit de montrer que si un élément $λ∈K$ appartient à 
$℘(\chap{K})+\chap{K}_{>0}$, il appartient également à
$℘(K)+K_{>0}$. Écrivons
$$
λ=\big(\frac{α₀}{β₀}-(\frac{α₀}{β₀})^p\big)+\sum_{θ≠0}
\big(\frac{α_θ}{β_θ}\big)^p b^θ,
$$
où les $α_\star$ sont dans $\chap{A}$, les $β_\star$ dans
$\chap{A}-\{0\}$, et $θ$ parcourt, comme dans \emph{loc. cit.}, l'ensemble
des applications $[1,r]→[0,p-1]$.
De façon équivalente, l'équation à coefficients dans $A$ 
$$
\big(\prod_{θ} Y_θ\big)^p λ = \big(\prod_{θ≠0} Y_θ \big)^p \cdot
\big( Y₀^{p-1} X₀-X₀^p\big) + \sum_{θ≠0} X_θ^p \big(\prod_{θ'≠θ} Y_{θ'}^p\big)
b^θ
$$
a pour solution 
$$
X_θ=α_θ, Y_θ=β_θ,
$$
où $θ$ parcourt cette fois-ci $[0,p-1]^{[1,r]}\cup \{0\}$.

D'après le théorème d'approximation ci-dessus, cette équation
a également une solution dans $A$, dont les coordonnées $Y$ sont 
non nulles, de sorte que $λ$ appartient bien à $℘(K)+K_{>0}$.
\end{démo}

\begin{rmr2}
Du fait que l'on peut définir un opérateur de Cartier
inverse $C^{-1}:Ω^i_A → Ω^i_{A}/\mathrm{d}Ω^{i-1}_A$ 
pour tout anneau $A$ de caractéristique $p$
(\cite{Nilpotent@Katz}, 7.2), on peut également déduire le résultat précédent
(en tout degré) de l'énoncé \ref{foncteur pf}. 
\end{rmr2}

\begin{crl2}\label{majoration complétion}
Sous les hypothèses de \ref{p-rang complété excellent}, on a l'inégalité 
$$
\dim_p(K)≤\dim_p(\chap{K}),
$$
où $p$ est la caractéristique du corps $K$.
\end{crl2}

\begin{démo}
Si $[K:K^p]<[\chap{K}:\chap{K}^p]$, il n'y a rien à démontrer
de sorte que l'on peut supposer les $p$-rangs finis, égaux à $r∈\NN$.
Il faut montrer que si $\HH^{r+1}_p(K')≠0$ pour une extension
finie $K'$ de $K$, il existe une extension finie $L/\chap{K}$
telle que $\HH^{r+1}_p(L)≠0$. D'après la proposition
précédente, il suffit de considérer le corps des fractions
du complété du normalisé de $A$ dans $K'$.
\end{démo}

(D'après \ref{énoncé}, cette inégalité est en fait une égalité.)

Terminons par une propriété d'invariance, élémentaire mais cruciale, 
de la $p$-dimension.

\begin{lmm2}\label{invariance p-dim}
Soient $k$ un corps, $k'/k$ une extension finie et $p$ un nombre premier.
Si $\dim_p(k)$ est \emph{fini}, on a égalité :
$$
\dim_p(k)=\dim_p(k').
$$
\end{lmm2}

\begin{démo}
Si $p$ est inversible sur $k$, c'est \cite{CG@Serre}, chapitre \textsc{ii},
¶4.1, proposition 10.
Si $p=\car{k}$, on sait déjà (\ref{invariance p-rang}) que les $p$-rangs
sont égaux. Il en résulte que l'on a une inégalité :
$\dim_p(k)≥\dim_p(k')$. Supposons $r=\prang(k)$ \emph{fini}.
Quitte à remplacer $k$ par une extension finie et $k'$ par une extension
composée, on est ramené à montrer que si $\HH^{r+1}_p(k')=0$,
on a également $\HH^{r+1}_p(k)=0$. Cela résulte de la surjectivité de la 
trace (\ref{rem trace}).
\end{démo}

\subsection{Le théorème de Katô}

\begin{thm2}[Kazuya Katô, \cite{Galois@Kato}]\label{kato}
Soit $A$ un anneau de valuation discrète hensélien excellent
de corps résiduel $k$ de caractéristique $p>0$ et de corps des
fractions $K$.
On a l'égalité 
$$
\dim_p(K)=1+\dim_p(k).
$$
\end{thm2}

Rappelons que l'hypothèse d'excellence signifie ici que l'extension
$\Frac{\,\chap{A}}/\Frac{\,A}$ est séparable.

Le lecteur pourra consulter avec profit l'exposé de J.-L. Colliot-Thélène \cite{Kato@Colliot-Thelene} pour une démonstration
du théorème ci-dessus.

Voici maintenant un corollaire du théorème précédent, qui
est l'analogue du (corollaire au) théorème de 
Tsen (\cite{CG@Serre}, chapitre \textsc{ii}, §3.3 et §4.2).

\begin{crl2}[Kazuya Katô et Takako Kuzumaki, \cite{Kato-Kuzumaki}]\label{cor kato}
Soient $K/k$ une extension de degré de transcendance $N$, et $p$ un nombre
premier.
On a l'inégalité
$$
\dim_p(K)≤N+\dim_p(k).
$$
\end{crl2}

(La démonstration du corollaire se fait par réduction 
au cas bien connu de la caractéristique nulle.)

\begin{rmr2}
Cette même méthode permet de déduire l'invariance de
la $p$-dimension par extension finie du théorème
de Katô.\end{rmr2}

Terminons cette section par une variante, élémentaire mais
plus explicite, du théorème \ref{kato}, qui semble laissée en exercice 
au lecteur dans la littérature (cf. p. ex.
\cite{Galois@Kato}). 

\begin{prp2}\label{Hr séries formelles}
Soit $k$ un corps de caractéristique $p>0$ et de $p$-rang $r<+∞$. Le
morphisme $$\HH^{r+1}_p(k)→\HH^{r+2}_p \big(k((t))\big)$$
envoyant la classe de $ω$ sur la classe de $ω\wedge \dlog(t)$
est un isomorphisme.
\end{prp2}

\begin{démo}
Soit $\{b₁,\dots,b_r\}$ une $p$-base du corps $k$ ; 
il résulte de \ref{p-rang séries formelles} et \ref{p-rang fractions}
que $\{b₁,\dots,b_r,t\}$ est une $p$-base de $k((t))$.
Considérons $k_{>0}$ et $k((t))_{>0}$ comme en \ref{structure Hr}
relativement à ces $p$-bases. 
(Rappelons que $k((t))_{>0}$ est un sous $k((t))^p$-espace
vectoriel de $k((t))$.)
Le morphisme $ω\mapsto ω\wedge \dlog(t)$
de l'énoncé correspond au morphisme canonique
$$
k/\Big(℘(k)+k_{>0}\Big)\sr{(\star)}{→} k((t))/\Big(℘\big(k((t))\big)+k((t))_{>0}\Big)
$$
(où $℘$ est le morphisme d'Artin-Schreier usuel)
déduit de l'inclusion $k\hra k((t))$ par passage au quotient.

Vérifions que $(\star)$ est une surjection. Puisque tout élément
de $k((t))$ de valuation strictement positive est dans l'image de 
$℘$ (cf. lemme \ref{surjectivité P} ci-dessus), il suffit de montrer
que tout élément $a_{-n}t^{-n}+\cdots a_{-1}t^{-1}+a_0$ est
congru modulo $℘\big(k((t))\big)+k((t))_{>0}$ à un élément de $k\subset k((t))$.
On peut supposer $a₀=0$.
Montrons par récurrence sur $n≥1$ que pour chaque $a∈k$, l'élément
$a t^{-n}$ appartient à $℘\big(k((t))\big)+k((t))_{>0}$. 
Par construction, pour tout $a∈k$, tout $r∈\ZZ$ et tout $i∈[1,p-1]$, l'élément
$(a t^{pr})\cdot t^{i}$ appartient à $k((t))_{>0}$ de sorte que
le résultat est acquis pour $n$ premier à $p$. 
Considérons maintenant le cas où $n$ est un multiple de $p$, $n=rp$.
Écrivons $a=a₀^p+a_{>0}$ où $a₀∈k$ et $a_{>0}∈k_{>0}$. 
On peut alors décomposer $a t^{-rp}$ en :
$$
\frac{a}{t^{rp}}=\big(\frac{a₀}{t^r}\big)^p + \frac{a_{>0}}{t^{rp}}.
$$
Le second terme, $\frac{a_{>0}}{t^{rp}}$, appartient à $k((t))_{>0}$ ; 
le premier terme est égal à $\frac{a₀}{t^r}-℘\big(\frac{a₀}{t^r}\big)$.
L'entier $r$ étant strictement inférieur à $rp$, on peut déduire de l'hypothèse
de récurrence que $\big(\frac{a₀}{t^r}\big)^p$ appartient à 
$℘\big(k((t))\big)+k((t))_{>0}$.

Vérifions maintenant que $(\star)$ est une injection.
Soit $a∈k$ tel que $a∈℘\big(k((t))\big)+k((t))_{>0}$
et montrons qu'il appartient à $℘(k)+k_{>0}$.
Écrivons
$a=℘(b^-+b₀+b^+)+b_{>0}$, où $b₀∈k$, $b^-$ (resp. $b^+$) est
un polynôme en $t^{-1}$ (resp. une série en $t$) sans terme constant, et $b_{>0}∈
k((t))_{>0}$. Puisque $℘(b^-)$ (resp. $℘(b^+)$)
est également un polynôme en $t^{-1}$ (resp. une série en $t$)
, sans terme constant, et que $℘(b₀)$ appartient à $℘(k)$, il suffit de vérifier que 
le terme constant de $b_{>0}$ appartient à $k_{>0}$.
Rappelons à cette fin que l'on a la décomposition :
$$
k((t))_{>0}=\Big(\bigoplus_{i∈[1,p-1],θ} k((t))^p b^θ t^i \Big)
\oplus \Big(\bigoplus_{θ≠0} k((t))^p b^θ\Big),
$$
où $θ$ parcourt l'ensemble des fonctions $[1,r]→[0,p-1]$.
La première somme directe ne contribue pas au terme constant donc
il suffit de vérifier que le terme constant d'un élément
de $\bigoplus_{θ≠0} k((t))^p b^θ$ appartient à $k_{>0}$ ; c'est clair.

\end{démo}

\section{Minoration de $\dim_p(K)$}\label{section minoration} 

Dans cette section, on démontre l'inégalité suivante, où $A$, $K$ 
et $k$ sont comme dans \ref{énoncé} :

\begin{eqnarray}\label{minoration}
\dim_p(K)≥\dim(A)+\dim_p(k).
\end{eqnarray}


\subsection{}Soit $A^ν$ le normalisé de $A$ dans son corps des fractions.
L'anneau $A$ étant excellent (cf. \ega{$\textsc{iv}_2$}{7.8}), il est également universellement japonais, 
de sorte que  $A^ν$ est fini sur $A$. 
Ce dernier étant hensélien et $A^ν$ étant intègre, l'anneau 
$A^ν$ est local. Il est également hensélien et excellent.

Soit $k^ν$ le corps résiduel de $A^ν$ ; c'est une extension finie
de $k$ de sorte que d'après \ref{invariance p-dim}, il 
suffit de démontrer \ref{minoration} dans le cas particulier
où $A$ est \emph{normal}. 

\subsection{Le cas d'inégale caractéristique (cf. \cite{Arithmetic@Saito}, §5)}

Soit $A$ comme ci-dessus, de dimension $≥2$, et 
$\got{p}$ un idéal premier de hauteur un. Soit $B$ le localisé
$A_{\got{p}}$ --- de corps des fractions $K$ ---, 
$\chap{B}$ son complété, et $L$ le corps des fractions
de $\chap{B}$. Puisque toute $L$-algèbre étale est induite
par une $K$-algèbre étale (\sga{4}{x}{2.2.1}), le morphisme
$G_L→G_K$ est une \emph{injection}.
Compte tenu de la décroissance de la dimension cohomologique
par passage à un sous-groupe fermé (\cite{CG@Serre}, chapitre I, §3.3,
proposition 14), on a donc :
$$
(\star\star) \ \cd_p(K)≥\cd_p(L).
$$
Le corps $L$ étant le corps des fractions d'un anneau de valuation 
discrète complet (donc hensélien, excellent), on a de plus l'inégalité
$$
\cd_p(L)≥1+\dim_p(κ(B)),
$$
où $κ(B)$ est le corps résiduel de $B$. Cela résulte du théorème
\ref{kato}. On achève la démonstration, par récurrence, en remarquant
que $κ(B)$ est le corps des fractions de l'anneau local intègre hensélien
excellent $A/\got{p}$, de dimension $\dim(A)-1$.

\subsection{Le cas d'égale caractéristique}

Procédant comme ci-dessus, il nous suffit de démontrer l'analogue
de $(\star\star)$ pour la $p$-dimension. C'est le contenu du lemme ci-dessous.

\begin{lmm2}
Soient $K$ un corps discrètement valué de caractéristique $p>0$ et $\chap{K}$
son complété. L'inégalité suivante est satisfaite :
$$
\dim_p(K)≥\dim_p(\chap{K}).
$$
\end{lmm2}

\begin{démo}
Si $K$ et $\chap{K}$ n'ont pas le même $p$-rang, il n'y a rien
à démontrer (cf. \ref{invariance p-rang}) ; supposons donc
qu'ils sont égaux à un entier $r<+∞$ et qu'il existe une extension finie $L'/\chap{K}$
telle que $\HH^{r+1}_p(L')≠0$. Il nous faut alors montrer qu'il
existe une extension finie $L/K$ telle que $\HH^{r+1}_p(L)≠0$.
Soit $k'$ le corps résiduel de $L'$ ; c'est une extension
finie du corps résiduel $k$ des corps discrètement valués 
$K$ et $\chap{K}$.
Il existe donc une extension finie $L/K$ telle
que $\OO_L/\MM_{\OO_K}=k'$ : on se ramène immédiatement
au cas où $k'/k$ est monogène, qui est bien connu
(cf. p. ex. \cite{CL@Serre}, chapitre I, proposition 15 ;
voir aussi \ega{$0_{\textsc{iii}}$}{10.3.1-2}).
En particulier, le complété $\chap{L}$ de $L$ est (abstraitement) isomorphe
à $L'$ (tous deux isomorphes au corps des séries de Laurent 
$k'((t))$).
D'après le lemme \ref{Hr complété}, pour un tel corps $L$,
le morphisme $\HH^{r+1}_p(L)→\HH^{r+1}_p(\chap{L})\isononcan
\HH^{r+1}_p(L')≠0$ est une
surjection. Ainsi, $\HH^{r+1}_p(L)≠0$, comme escompté.
\end{démo}

\section{Majoration de $\dim_p(K)$ : le cas d'égale 
caractéristique}\label{section majoration egale} 

Dans cette section, on démontre par récurrence sur la dimension 
l'inégalité suivante, où $A$, $K$ 
et $k$ sont comme dans \ref{énoncé}, et où l'on suppose de plus \emph{$A$ de
caractéristique $p>0$} :

\begin{eqnarray}\label{maj égale}
\dim_p(K)≤\dim(A)+\dim_p(k).
\end{eqnarray}

Il résulte de \ref{minoration} et  \ref{majoration complétion} que l'on peut supposer 
$A$ normal complet. (On utilise également le fait que 
le complété d'un anneau local excellent normal est normal (\ega{iv}{7.6.1}).)

Notons $d$ la dimension de $A$, $k$ son corps résiduel, $r$ le $p$-rang
de $k$, que l'on peut supposer fini, et $K$ le corps des fractions de $A$.
Commençons par un lemme élémentaire.

\begin{lmm}
Sous les hypothèses précédentes, on a :
$$
\prang(K)=d+r.
$$
\end{lmm}

\begin{démo}
D'après le théorème de structure de Cohen \ega{$0_{\textsc{iv}}$}{19.8.8 ii},
il existe un sous-anneau $A₀$ de $A$, isomorphe à $k[[t₁,\dots,t_d]]$,
tel que le morphisme $\SP(A)→\SP(A₀)$ soit \emph{fini}
de sorte que $\prang\,K=\prang\,\Frac(A₀)$ (\ref{invariance p-rang}).
D'après \ref{p-rang séries formelles}, le terme de droite est égal à $d+r$.
\end{démo}

Il nous faut donc montrer que si $\dim_p(k)=r$ (\cad
si $\HH^{r+1}_p(k')=0$ pour toute extension finie $k'/k$),
hypothèse que nous allons maintenant supposer satisfaite,
on a également $\HH^{d+r+1}_p(K')=0$ pour toute extension finie
$K'/K$. 
Il suffit de montrer que $\HH^{d+r+1}_p(K)=0$.

\subsection{}\label{omega}Posons $n=d+r$ et considérons un élément de $Ω^n_K$, 
que l'on écrit $\frac{ω}{f}$, où $ω∈Ω^n_A/\textrm{torsion}$ et 
$f∈\MM_A-\{0\}$. 
D'après le théorème de structure de Cohen-Gabber 
(\ref{gabber}), il existe un sous-anneau $A₀$ de
$A$, isomorphe à $k[[x₁,\dots,x_d]]$ tel que le morphisme
$π:\SP(A)→\SP(A₀)$ soit fini et \emph{génériquement étale}.
Soit $f₀=\mathrm{N}_{X/X₀}(f)∈A₀$ la norme de l'élément $f$ 
(\ega{ii}{6.5}). L'élément $f₀$ divise $f$ de sorte
que l'on peut supposer, et l'on supposera, que l'élément $f$ appartient à $A₀$.
Considérons le fermé $R$ de $X₀:=\SP(A₀)$ au-dessus duquel le morphisme
$π$ est ramifié. Il existe un élément non nul $a∈A₀$ tel que $R$ soit contenu
dans le fermé $V(a)$. Enfin, quitte à remplacer $a$ et $f$ par $af$, on peut supposer
$a=f$.

Rappelons le théorème de préparation, qui nous permettra de 
rendre le lieu de ramification fini, par projection, sur un schéma
de dimension un de moins.

\begin{thm}[Théorème de préparation de Weierstraß, \ac{vii}{3}{7-8}{}]\label{weierstrass} \ \\ 
Soient $A$ un anneau local séparé complet d'idéal maximal $\MM$
et $n$ un entier.
\begin{enumerate}
\item Soit $k$ un entier et $f∈A[[\sous{X},T]]$
($\sous{X}=\{X₁,\dots,X_n\}$) une série entière
\emph{$k$-régulière} relativement à $T$, \cad congrue à $(u∈A[[T]]^×)\cdot
T^k$ modulo $(\MM,\sous{X})$. 
Pour tout $g∈A[[\sous{X},T]]$, il existe un unique
couple $(q,r)∈A[[\sous{X},T]]×A[[\sous{X}]][T]$ tel que
$g=qf+r$ et $\deg_T(r)<k$.
\item Soit $k$ un entier. Si $f∈A[[\sous{X},T]]$ est $k$-régulière
relativement à $T$, il existe un polynôme
$P=T^k+∑_{i<k} p_i T^i$, où
$p_i∈(\MM,\sous{X})A[[\sous{X}]]$, et une unité
$u∈A[[\sous{X},T]]^×$ tels que $f=uP$.

\item Soit $f∈A[[\sous{X},T]]$ non nulle modulo $\MM$.
Il existe un entier $k$ et un automorphisme
$A[[T]]$-linéaire $c$ de $A[[\sous{X},T]]$, tel que
$c(X_i)=X_i+T^N_i$ et $c(f)$ soit $k$-régulier.

\end{enumerate}
\end{thm}

D'après (iii) et (ii), quitte à changer de coordonnées, on peut supposer 
$a$ et $f$ égaux à un même polynôme \emph{unitaire} en $x_d$, 
dont les autres coefficients appartiennent à l'idéal
maximal de $k[[x₁,\dots,x_{d-1}]]$. 
On notera $G$ ce polynôme, qui
appartient en particulier à l'anneau $\gtilde{A₀}:=k[[x₁,\dots,x_{d-1}]]\{x_d\}$,
hensélisé de $k[[x₁,\dots,x_{d-1}]][x_d]$ en l'origine.

\begin{lmm}\label{lemme complétion}
Soient $B$ un anneau local complet noethérien et $P∈B[X]$ un polynôme
de la forme $X^n+∑_{i<n}b_i X^i$, où $b_i∈\MM_{B}$.

\begin{enumerate}
\item Le complété $(P)$-adique de $B\{X\}$ s'identifie à
$B[[X]]$.
\item La paire $\big(\SP(B\{X\}),V(P)\big)$ est \emph{hensélienne}.
\end{enumerate}
\end{lmm}

\begin{démo}
Vérifions (i). Soit $Q∈
B[X]$ un polynôme unitaire. Il résulte de \ref{weierstrass} (i), que l'anneau quotient
$B[[X]]/(Q)$ est 
isomorphe comme $B$-module à $B[X]/(X^{\deg(Q)})$
et en particulier fini sur $B$. Son sous-anneau (par
fidèle platitude)
 $B\{X\}/(Q)$ est
donc également fini sur $B$, donc complet, et finalement isomorphe à $B[[X]]/(Q)$.
Utilisant ce fait pour $Q=P^n$, $n≥0$, on en déduit
que le séparé-complété $(P)$-adique de $B\{X\}$ est 
isomorphe 
à celui de $B[[X]]$ ; 
ce dernier est isomorphe à $B[[X]]$ si l'on suppose de plus que
$P$ n'est pas constant.

Vérifions (ii).
Rappelons qu'une paire $(\SP(C),V(I))$ est dite \emph{hensélienne}
si pour tout polynôme $f∈C[T]$, toute racine \emph{simple} de $g$ dans $C/I$
se relève en une racine dans $C$. On laisse le soin au lecteur de vérifier
que pour tout anneau local hensélien $C$ et tout idéal $I\subset \MM_C$,
la paire $(\SP(C),V(I))$ est hensélienne. On applique alors ce résultat
à $C=B\{X\}$ et $I=(P)$.
\end{démo}
Terminons ces rappels par l'énoncé du théorème d'algébrisation suivant :

\begin{thm}[Renée Elkik, \cite{Solutions@Elkik}, théorème 5]\label{elkik}
Soient $(X=\SP(A),Y=V(I))$ une paire hensélienne avec $A$ noethérien,
et $U$ le sous-schéma ouvert complémentaire de $Y$ dans $X$.
Notons $X_{\chap{Y}}$ le complété de $X$ le long de $Y$, $\chap{Y}$
le fermé correspondant à $Y$ et $\chap{U}$ son complémentaire dans 
$X_{\chap{Y}}$.
Le foncteur $X'↦X'×_X X_{\chap{Y}}$ induit une
équivalence de catégories entre la catégorie des
$X$-schémas finis, étales sur $U$, et la catégorie
des $X_{\chap{Y}}$-schémas finis, étales sur $\chap{U}$. 
\end{thm}

Reprenons les notations en vigueur après le théorème \ref{weierstrass}.
Le fermé $V(G)$, contenant le lieu de ramification, étant 
défini par une équation polynomiale unitaire à coefficients dans
$k[[x₁,\dots,x_{d-1}]]$, il résulte du lemme \ref{lemme complétion}
et du théorème \ref{elkik} que l'on peut algébriser le morphisme
$π:\SP(A)→\SP(A₀)$. En d'autres termes, on a un diagramme cartésien :
$$\xymatrix{
\ar @{} [dr] |{\square}
X=\SP(A) \ar[d]^{π} \ar[r] & \gtilde{X}=\SP(\gtilde{A}) \ar[d]^{\text{fini, gén.
étale}} \\
X₀=\SP(A₀) \ar[r] & \gtilde{X₀}=\SP(\gtilde{A₀})
}
$$
où, rappelons-le, $\gtilde{A₀}$ est $k[[x_1,\dots,x_{d-1}]]\{x_d\}$.

D'après le lemme \ref{lemme complétion}, appliqué 
à $B=k[[x₁,\dots,x_{d-1}]]$ et $P=G$, l'anneau complet $A$, fini sur $A₀$,
s'identifie au complété $(G)$-adique de $\gtilde{A}$.
On peut donc décomposer
$\frac{ω}{G}$ en 
$$
\frac{ω}{G}=\frac{\gtilde{ω}}{G}+ω'
$$
où $\gtilde{ω}∈Ω^n_{\gtilde{A}}/\mathrm{torsion}$ et 
$ω'∈GΩ^n_A/\mathrm{torsion}\subset \MM_A Ω^n_A/\mathrm{torsion}$.

Le corps des fractions $\gtilde{K}$ de $\gtilde{A}$ est de degré de transcendance
un sur le corps des fractions $L$ de $k[[x₁,\dots,x_{d-1}]]$.
Par hypothèse de récurrence sur la dimension,
et compte tenu du fait que l'on a supposé
$\dim_p(k)=\prang(k)$,
on a $$\dim_p(L)≤(d-1)+\dim_p(k)=d-1+r.$$ 
Enfin, d'après le corollaire au théorème de Katô (\ref{cor kato}),
on a $$\dim_p(\gtilde{K})≤\dim_p(L)+1.$$
Finalement, $\dim_p(\gtilde{K})≤d+r$.
Ainsi, la classe de $\frac{\gtilde{ω}}{G}$ dans $\HH^{r+1}_p(\gtilde{K})$
est nulle ; \emph{a fortiori} son image dans $\HH^{r+1}_p(K)$ l'est également.

On achève la démonstration en remarquant que la classe $ω'$,
$\MM$-adiquement proche de zéro, est elle-aussi nulle dans $\HH^{r+1}_p(K)$ :
cela résulte d'une variante de \ref{surjectivité P}, laissée
en exercice au lecteur.

\section{Majoration de $\dim_p(K)$ : le cas d'inégale
caractéristique}\label{section majoration mixte}

\subsection{}\label{notations mixte}Dans cette section, on démontre 
par récurrence sur la dimension l'inégalité suivante, où $A$, $K$ 
et $k$ sont comme dans \ref{énoncé}, et où l'on suppose de plus $A$ de
caractéristique nulle :

\begin{eqnarray}\label{majoration mixte}
\mathrm{cd}_p(K)≤\dim(A)+\dim_p(k).
\end{eqnarray}

Il résulte de \ref{minoration} que l'on peut supposer 
$A$ normal. Comme en \ref{majoration complétion}
on utilise le théorème d'approximation de Popescu, sous la forme suivante, 
pour se ramener au cas complet.

\subsection{Réduction au cas complet}

\begin{lmm2}\label{lemme gabber}
Soit $A$ un anneau local hensélien excellent intègre de corps des fractions $K$.
Soient $\chap{A}$ son complété et $\chap{K}$ son corps des fractions.
Il existe un ensemble filtrant de sous-$K$-algèbres de type fini 
$(B_i)_{i∈I}$ de $\chap{K}$ tel que $\chap{K}=\colim \,B_i$ et chaque morphisme
$\SP(B_i)→\SP(K)$ ait une \emph{section}. 
\end{lmm2}

\begin{crl2}\label{foncteur pf}
Soient $A$ comme en \ref{lemme gabber} et $\mathsf{F}$ un foncteur de localement présentation finie
$(A-\textrm{algèbres})→\Ens$.
Le morphisme canonique $\mathsf{F}(K)→\mathsf{F}(\chap{K})$
est une \emph{injection}.
\end{crl2}

D'après \sga{4}{vii}{5.7}, on en déduit :

\begin{crl2}
Soient $A$ comme en \ref{lemme gabber}, $p$ un nombre premier
et $i$ un entier.
Le morphisme canonique
$$
\HH^{i}(\SP(K)_{\et},\ZZ/p)→\HH^{i}(\SP(\chap{K})_{\et},\ZZ/p)
$$
est une injection.
\end{crl2}

\begin{démo}[Démonstration du lemme \ref{lemme gabber}]
Il suffit de démontrer que pour toute sous-$K$-algèbre de type fini
$B$ de $\chap{K}$, le morphisme $\SP(B)→\SP(K)$ a une section.
L'inclusion $K$-linéaire $B\ra \chap{K}$ correspond à un $\chap{K}$-point
d'un système d'équations convenables $f_1,\dots,f_r$ (définissant
$B$ sur $K$) que l'on peut supposer à coefficients dans $A$. Quitte 
à chasser comme en \ref{Hr complété excellent} (démonstration) 
les dénominateurs, on peut appliquer \ref{popescu} pour obtenir un $K$-point,
c'est-à-dire une section du morphisme $\SP(B)→\SP(K)$.
\end{démo}

\subsection{}\label{Epp et maj grossière} En plus des hypothèses en 
vigueur dans \ref{notations mixte}, on suppose maintenant $A$ \emph{complet}. 

Pour pouvoir utiliser le théorème de structure de Cohen-Gabber (\ref{gabber})
--- qui suppose la « fibre spéciale » réduite ---,
nous ferons usage du théorème
suivant :

\begin{thm2}[Helmut Epp, \cite{Eliminating@Epp},
théorème 1.9]\label{Epp}
Soit $T→S$ un morphisme local dominant de traits complets, de caractéristique
résiduelle $p>0$. Notons $κ_S$ et $κ_T$ leurs corps résiduels respectifs.
Supposons $κ_S$ \emph{parfait} et le sous-corps parfait maximal
de $κ_T$ \emph{algébrique} sur $κ_S$.
Il existe une extension finie de traits $S'→S$ telle que
le produit fibré réduit normalisé
$$
T':=(T\times_S S')_{\red}^{ν}
$$
ait une \emph{fibre spéciale réduite} au-dessus de $S'$.
\end{thm2}

\begin{rmr2}\label{rmr Epp}
En caractéristique mixte, on vérifie immédiatement que le produit fibré $T\times_S S'$ 
est réduit 
et que si l'on suppose seulement $S$ complet (mais pas
nécessairement $T$), la conclusion est encore valable.
(Cf. \emph{loc. cit.}, §2 pour le cas général.)
\end{rmr2}

\subsubsection{}\label{stabilité E}
Commençons par vérifier que l'hypothèse sur les corps résiduels
est satisfaite dans de nombreux cas. Nous dirons qu'une extension
de corps $K/k$ de caractéristique $p>0$ a la \emph{propriété de Epp} si tout
élément du sous-corps $K^{p^∞}:=\cap_{i≥0} K^{p^{n}}$ de $K$ 
est algébrique séparable sur $k$.

\begin{lmm3}
Pour tout corps $K$ de caractéristique $p>0$, on a, dans une clôture séparable $K^{\sep}$ de $K$, 
$$(K^{p^∞})^{\sep}=(K^{\sep})^{p^∞}.$$
\end{lmm3}

\begin{démo}
L'inclusion $(K^{p^{∞}})^{\sep}\subset (K^{\sep})^{p^∞}$ est 
évidente : $K^{p^∞}$ est parfait donc toute extension algébrique,
en particulier sa clôture séparable
$(K^{p^{∞}})^{\sep}$, l'est également. Comme cette dernière est contenue
dans $K^{\sep}$, elle est également contenue dans son plus grand
sous-corps parfait $(K^{\sep})^{p^∞}$.

Réciproquement, considérons $x∈(K^{\sep})^{p^∞}$, et notons, pour chaque
entier $n≥0$, $x_n$ sa racine $p^n$-ième dans $K^{\sep}$ et $f_n$ son polynôme
minimal (\emph{unitaire}). Compte tenu d'une part de l'expression de $f_n$ en fonction
des polynômes symétriques en les conjugués galoisiens de $x_n$ et d'autre part
de l'injectivité et de l'additivité de l'élévation à la puissance $p^n$-ième,
on a l'égalité $f₀=f_n^{(p^n)}$, où $f_n^{(p^n)}$ est le polynôme
obtenu à partir de $f_n$ en élevant les coefficients à la puissance $p^n$-ième.
Il en résulte que les coefficients du polynôme minimal $f₀$ de $x$
appartiennent à $K^{p^∞}$.

\end{démo}

\begin{prop3}[Cf. \cite{Eliminating@Epp}, §0.4]\label{prop Epp}
\ \\ 
\begin{enumerate}
\item Soient $L/K$ et $K/k$ ayant la propriété
de Epp. Alors, $L/k$ a la propriété de Epp.
\item Toute extension finie a la propriété de Epp.
\item Pour tout entier $d$, et tout corps
$k$, l'extension $\big(\Frac{\,k[[x₁,\dots,x_d]]}\big)/k$
a la propriété de Epp.
\item Soient $A$ un anneau local complet noethérien intègre,
et $k\subset A$ un corps de représentants (s'envoyant 
isomorphiquement sur le corps résiduel de $A$). Alors, l'extension $(\Frac{\,A})/k$
a la propriété de Epp.
\end{enumerate}
\end{prop3}

\begin{démo}
(i) Par hypothèse, $L^{p^∞}\subset K^{\sep}$. Comme
le corps $L^{p^∞}$ est \emph{parfait}, on en déduit 
que $L^{p^∞}\subset (K^{\sep})^{p^∞}=(K^{p^∞})^{\sep}\subset
k^{\sep}$, où l'égalité résulte du lemme précédent. 

(ii) Toute extension étale a tautologiquement
la propriété de Epp. D'après (i), il reste à considérer
le cas d'une extension radicielle $K/k$.
Si elle est de hauteur $≤r$, on a $K^{p^r}\subset k$ et
en particulier $K^{p^∞}\subset k\subset k^{\sep}$.

(iii) Soit $A=k[[x₁,\dots,x_d]]$ et $K$ son corps 
des fractions. Montrons que $K^{p^∞}=k^{p^∞}$.
Comme $K$ est contenu dans $k((x₁,\dots,x_{d-1}))((x_d))$,
on se ramène par récurrence au cas où $d=1$.
Tout élément non nul de $k((t))^{p^∞}$ a une valuation infiniment
$p$-divisible donc nulle, de sorte que
$k((t))^{p^∞}-\{0\}$ est contenu dans 
$k[[t]]^×$ et finalement dans $k^{p^∞}$ par un calcul immédiat.

(iv) Cela résulte des observations précédentes et du théorème
de structure de Cohen.

\end{démo}

\subsubsection{}Soit $k₀=k^{p^∞}$ le sous-corps
parfait maximal du corps résiduel $k$ de $A$ et notons $W₀=W(k₀)$
l'anneau des vecteurs de Witt correspondant.
Il résulte du théorème de Cohen que l'on a un morphisme
$X:=\SP(A)→S₀:=\SP(W₀)$ relevant l'inclusion $k₀\inj k$
(\ega{$0_{\textsc{iv}}$}{19.8.6}).

Pour tout point maximal $\got{p}$ de la fibre spéciale $X_p$,
l'anneau de valuation discrète $A_{\got{p}}$ a pour corps
résiduel $\Frac{\,A/\got{p}}$, où $A/\got{p}$ est un anneau complet
intègre noethérien de corps résiduel $k$. D'après
\ref{prop Epp} (i) \& (iv), l'extension $\Frac(A/\got{p})/k₀$ satisfait donc
l'hypothèse du théorème de Epp. Les idéaux $\got{p}$ étant
en nombre fini et la conclusion du théorème de Epp
(\ref{Epp} et \ref{rmr Epp}) étant
stable par changement de base fini (c'est un résultat
de lissité formelle), il existe donc un changement de base fini $S₀'=\SP(W₀')→S₀$ 
tel que la fibre spéciale du produit fibré normalisé $X':=(X\times_{S₀}
S₀')^ν=\SP(A')$
soit réduite en ses points maximaux. (On utilise le fait que les points maximaux
de la fibre spéciale de $X'→ S₀'$ se trouvent au-dessus des points maximaux de la fibre
spéciale de $X→S₀$ ; cf. p. ex. \ega{$0_{Ⅳ}$}{16.1.6}.)

D'après le lemme suivant, la fibre spéciale est alors 
réduite. 

\begin{lmm3}
Soit $X$ un schéma noethérien normal. Tout diviseur de Cartier effectif 
génériquement réduit est réduit.
\end{lmm3}

\subsubsection{}\label{fibration mixte}
Notons $k₀'$ le corps résiduel de $W₀'$, fini
sur $k₀$, $\varpi'$ une uniformisante de $W₀'$, et considérons
une composante connexe $X''=\SP(A'')$ de $X'$. Soit $k''$ son
corps résiduel, fini sur $k$.
L'inclusion $k₀'\inj k''$ déduite du morphisme $X''→S'$
est formellement lisse, car $k₀'$ est parfait, donc
se relève d'après \ega{$0_{\textsc{iv}}$}{19.7. 1 et 2} en
un morphisme \emph{formellement lisse} $W'₀→I''$ où 
$I''$ est un anneau local complet noethérien. Cet anneau
est un anneau de valuation discrète.
L'anneau $A''/\varpi'$ étant réduit (et équidimensionel de dimension $d-1$ de corps
résiduel $k''$), il existe d'après le théorème de Cohen-Gabber (\ref{gabber}),
un relèvement ($k₀'$-linéaire) $k''\hra A''/\varpi'$ et 
des éléments $x_1,\dots,x_{d-1}$ dans l'idéal maximal de $A''/ϖ'$ 
tels que le morphisme induit
$k''[[t_1,\dots,t_{d-1}]]\ra A''/ϖ'$, envoyant
l'indéterminée $t_i$ sur $x_i$, soit
fini, \emph{génériquement étale} en haut et en bas.

Par lissité formelle de $W₀'→I''$, le morphisme composé
$I''→k''→A''/\varpi'$ se relève en un $W₀'$-morphisme $I''→A''$.
En relevant les $x_i$ dans $A''$,
cela nous permet de construire un morphisme
$A₀'':=I''[[t_1,\dots,t_{d-1}]]→A''$,
fini injectif (cf. p. ex. \ega{$0_{\textsc{iv}}$}{19.8.8 (démonstration)}), 
\emph{étale} au-dessus du point générique de la fibre spéciale.

Notons $K''$ le corps des fractions de $A''$.

\subsection{}\label{majoration cas corps} Le but de ce paragraphe est de démontrer la
proposition suivante :

\begin{prp2}\label{cas réduit} 
Les groupes $\HH^{N}(\SP(K'')_{\et},\ZZ/p)$ sont nuls pour $N>\dim(A)+\dim_p(k)$.
\end{prp2}

Remarquons que $\dim(A)=\dim(A'')$ et que $\dim_p(k)=\dim_p(k'')$.

Soit $N$ comme ci-dessus et considérons une classe $c∈\HH^N(\SP(K'')_{\et},\ZZ/p)$.
Nous allons commencer par montrer que $c$ s'étend à un grand ouvert :

\begin{lmm2}\label{gros ouvert}
Il existe un ouvert $U\subset X''$ 
\emph{contenant les points maximaux de la fibre spéciale} et
une classe $c_{U}∈\HH^N(U_{\et},\ZZ/p)$ s'envoyant sur $c$ par restriction
à $\SP(K'')$.
\end{lmm2}

\begin{démo}[Démonstration du lemme]
Soit $\got{p}$ un point maximal de la fibre spéciale
$\SP(A''/\varpi')$. L'anneau $A''$ étant normal,
le localisé $A''_{\got{p}}$ est un anneau de valuation discrète.
Soit ${K''}_{\got{p}}^h$ le corps des fractions de l'hensélisé de $A''_{\got{p}}$
et $c_{\got{p}}$ la restriction de $c$ à $\SP({K''}_{\got{p}}^h)$.
Il résulte du théorème de Katô (\ref{kato}) que l'on a 
$\mathrm{cd}_p({K''}_{\got{p}}^h)=1+\dim_p(\Frac{\,A''/\got{p}})$.
D'après \ref{maj égale}, 
on a donc $\mathrm{cd}_p({K''}_{\got{p}}^h)≤\dim(A'')+\dim_p(k'')<N$, de sorte
que $c_{\got{p}}=0$.
Joint au lemme ci-dessous, cela montre que la classe 
$c$ appartient à l'image du morphisme de restriction 
$\HH^{N}(\SP(A''_{\got{p}})_{\et},\ZZ/p)→\HH^{N}(\SP(K'')_{\et},\ZZ/p)$.
Il existe donc, pour chaque $\got{p}$ comme ci-dessus, 
un ouvert $U_{\got{p}}$ de $X''$ contenant $\got{p}$
et une classe $C_{\got{p}}$ sur $U_{\got{p}}$ induisant $c$ sur $\SP(K)$.
Quitte à rétrécir ces ouverts, on peut utiliser inductivement la suite
exacte de Mayer-Vietoris pour recoller ces $C_{\got{p}}$ en une classe
$C_U$ sur $U=\cup U_{\got{p}}$.

\end{démo}

\begin{lmm2}
Soit $A$ un anneau de valuation discrète et notons $K$ son corps des fractions.
Considérons $A^h$ son hensélisé, $K^h$ son corps des fractions et 
$c∈\HH^{N}(\SP(K)_{\et},\ZZ/n)$, où $(n,N)∈\NN^2$.
Si l'image de $c$ dans $\HH^N(\SP(K^h)_{\et},\ZZ/n)$ est nulle, la classe
$c$ appartient à l'image du morphisme de restriction
$\HH^N(\SP(A)_{\et},\ZZ/n)→\HH^N(\SP(K)_{\et},\ZZ/n)$.
\end{lmm2}

\begin{démo}
Le morphisme $\SP(A^h)→\SP(A)$ induisant un isomorphisme
sur les localisations strictes et un isomorphisme au-dessus 
du point fermé $s$ de $\SP(A)$, on voit facilement d'après
\sga{4}{v}{6.4 et 6.5} et \sga{4}{viii}{5.2} que 
le morphisme d'adjonction $\RG_{\{s\}}(\SP(A)_{\et},\ZZ/n)→
\RG_{\{s^h\}}(\SP(A^h)_{\et},\ZZ/n)$ (où $s^h$ est le point
fermé de $\SP(A^h)$) est un isomorphisme.
Le morphisme de suites exactes 
$$
\xymatrix{
\HH^{N}(\SP(A)_{\et},\ZZ/n)\ar[d]\ar[r] &
\HH^{N}(\SP(K)_{\et},\ZZ/n)\ar[d]\ar[r] & \HH^{N+1}_{\{s\}}(\SP(A)_{\et},\ZZ/n)
\ar[d]^{\mathrm{isom.}}\\
\HH^{N}(\SP(A^h)_{\et},\ZZ/n) \ar[r] &
\HH^{N}(\SP(K^h)_{\et},\ZZ/n) \ar[r] & \HH^{N+1}_{\{s^h\}}(\SP(A^h)_{\et},\ZZ/n) 
}
$$
permet alors de conclure.
\end{démo}

Revenons à la démonstration de \ref{cas réduit}. On
peut supposer $d≥2$, sans quoi le résultat est déjà
connu.
Notons $π:X''=\SP(A'')→X''₀=\SP(A''₀)$ le morphisme 
considéré en \ref{fibration mixte}.
Il est fini et \emph{étale} au-dessus du complémentaire d'un fermé $F_π$ de $X''₀$ 
\emph{ne contenant pas} la fibre spéciale (au-dessus de $S₀'$).
Notons $F_c$ le fermé $π(X''-U)\subset X''₀$, où $U$ est comme en \ref{gros ouvert},
et posons $F:=F_π\cup F_c$. Par hypothèse, $F$ est contenu dans le lieu 
d'annulation d'une fonction $f∈A''₀$ telle que $f\notin (\varpi')$.
Ainsi, d'après \ref{weierstrass} (ii) et (i) (appliqué
à $κ=I'$ et $\MM=(\varpi')$), on peut supposer
que $f$ appartient à $I'[[x₁,\dots,x_{d-2}]][x_{d-1}]$ et est \emph{unitaire}
en $x_{d-1}$. 

De même qu'en égale caractéristique, il résulte du lemme 
\ref{lemme complétion} et du théorème \ref{elkik} que l'on peut algébriser le morphisme
$π:\SP(A'')→\SP(A'₀)$ : il existe un diagramme cartésien
$$\xymatrix{
\ar @{} [dr] |{\square}
X''=\SP(A'') \ar[d]^{π} \ar[r] & \gtilde{X''}=\SP(\gtilde{A''}) \ar[d]^{\text{fini, gén.
étale}} \\
X''₀=\SP(A''₀) \ar[r] & \gtilde{X''₀}=\SP(\gtilde{A''₀})
}
$$
où $\gtilde{A''₀}$ est $I''[[x_1,\dots,x_{d-2}]]\{x_{d-1}\}$. 

La paire $(\gtilde{X''},V(f))$ est hensélienne, car $(\gtilde{X''₀},V(f))$
l'est, de sorte qu'il résulte du théorème de comparaison de Fujiwara-Gabber
(\cite{Tubular@Fujiwara}, 6.6.4) que 
le morphisme
$$
\HH_{\et}^N(\gtilde{X''}-V(f),\ZZ/p)→\HH^N_{\et}(X''-V(f),\ZZ/p)
$$
est un \emph{isomorphisme}.
Ainsi, la classe $c∈\HH^N(\SP(K'')_{\et},\ZZ/p)$, qui a été préalablement
étendue à $X''-V(f)$, provient, par restriction, d'un élément
de $\HH^N_{\et}(\gtilde{X''}-V(f),\ZZ/p)$. 
Soient $\gtilde{K''}$ le corps des fonctions rationnelles du schéma intègre
$\gtilde{X''}$ et $L$ le corps des fractions de l'anneau $I''[[x₁,\dots,x_{d-2}]]$. 
L'extension $\gtilde{K''}/L$ est de degré de transcendance un de sorte
que $\mathrm{cd}_p(\gtilde{K''})≤1+\mathrm{cd}_p(L)$ (\ref{cor kato}). D'après l'hypothèse
de récurrence, on sait d'autre part que
$\mathrm{cd}_p(L)≤(\dim(A'')-1)+\dim_p(k'')$.
Finalement, $\mathrm{cd}_p(\gtilde{K})≤\dim(A)+\dim_p(k)$ et
$c$, qui appartient à l'image de $\HH^{N}(\SP(\gtilde{K''})_{\et},\ZZ/p)→
\HH^N(\SP(K'')_{\et},\ZZ/p)$ est nul. (Rappelons que $N>\dim(A)+\dim_p(k)$.)

Ceci achève la démonstration de la proposition \ref{cas réduit}.

Soient $K₀$ (resp. $K'₀$) le corps des fonctions de $S₀$ (resp. $S'₀$).
Par construction, il résulte de \ref{cas réduit} que
les groupes de cohomologie $\HH^{N}(\SP(K\otimes_{K₀} K₀')_{\et},\ZZ/p)$ sont nuls pour
$N>\dim(A)+\dim_p(k)$.

\begin{crl2}\label{intermédiaire}
Soit $\sur{K₀}$ une clôture séparable de $K₀$. Les groupes
$\HH^{N}(\SP(K\otimes_{K₀}\sur{K₀})_{\et},\ZZ/p)$ sont nuls pour $N>\dim(A)+\dim_p(k)$.
\end{crl2}

\begin{démo}
Soit en effet $\gtilde{W}/W₀'$ une extension finie.
La fibre spéciale du morphisme $\SP(A'\otimes_{W₀'}
\gtilde{W})→\SP(\gtilde{W})$ est encore réduite, de sorte que l'on
peut appliquer le théorème de Epp et la proposition précédente
pour en déduire la nullité de 
$\HH^N(\SP(K\otimes_{K₀} \Frac{(\gtilde{W})})_{\et},\ZZ/p)$.
On passe alors à la limite.

\end{démo}

\begin{prp2}\label{maj grossière}
Soient $A$, $K$, $k$ et $p$ comme dans \ref{énoncé}. Supposons $A$ normal, complet
de caractéristique mixte.
On a alors l'inégalité :
$$
\mathrm{cd}_p(K)≤\dim(A)+\dim_p(k)+2.
$$
\end{prp2}

Nous utiliserons le lemme élémentaire suivant :

\begin{lmm2}\label{coeff constants}
Soient $K$ un corps, $N$ un entier et $p$ un nombre premier.
Supposons que pour toute extension finie \emph{séparable} $L$ de $K$
on ait $\HH^M(\SP(L)_{\et},\ZZ/p)=0$ pour tout $M>N$. Alors,
$\mathrm{cd}_p(K)≤N$.
\end{lmm2}

\begin{démo}[Démonstration du lemme]
Soient $\sur{K}$ une clôture séparable de $K$ et $H$ un $p$-Sylow de
$\mathrm{Gal}(\sur{K}/K)=:G_K$. On a l'égalité
$\mathrm{cd}_p(G_K)=\mathrm{cd}_p(H)$ (\cite{CG@Serre}, chap. I, 
§3.3, corollaire 1). Puisque que $H$ est un $p$-groupe,
$\mathrm{cd}_p(H)$ est le plus grand entier $D$ (s'il existe) tel que
$\HH^M(H,\ZZ/p)=0$ pour tout $M>D$ (\emph{loc. cit.}, §4, proposition 21). 
Or, pour chaque $M$,
le groupe $\HH^M(H,\ZZ/p)=\HH^M(\sur{K}^H,\ZZ/p)$ est une colimite
de groupes $\HH^M(L,\ZZ/p)$ avec $L/K$ finie étale. Sous nos hypothèses,
les termes de cette colimite sont tous nuls pour $M>N$ de sorte que 
$\mathrm{cd}_p(H)≤N$ comme escompté. 
\end{démo}

\begin{démo}[Démonstration de la proposition]
Rappelons que le corps $K₀$ est un corps local à corps résiduel \emph{parfait}
de sorte que $\mathrm{cd}_p(K₀)≤2$ (\emph{loc. cit.}, chap. II, 4.3). 
Cette observation, jointe
au résultat l'annulation \ref{intermédiaire}, permet de 
déduire de la suite spectrale 
$$
E^{i,j}_2=\HH^i(G_{K₀},\HH^j_{\et}(K\otimes_{K₀}\sur{K₀},\ZZ/p))⇒
\HH_{\et}^{i+j}(K,\ZZ/p)
$$
l'annulation des groupes 
$\HH_{\et}^{N}(K,\ZZ/p)$ pour $N>\dim(A)+\dim_p(k)+2$.
Ceci étant également valable pour les extensions finies de $K$, 
on peut utiliser le lemme ci-dessus pour conclure.
\end{démo}

\subsection{Fin de la démonstration}\label{fin corps} Fixons $d$ et $r$ deux entiers et considérons
le plus petit entier $N$ tel que 
pour tout $A$ comme en \ref{maj grossière}, de dimension $d$
et de corps résiduel de $p$-dimension $r$,
on ait $\mathrm{cd}_p(\Frac\, A)≤N$. 
On a vu ci-dessus qu'un tel entier $N$ existe (et est inférieur à 
$d+r+2$). D'après le lemme \ref{coeff constants},
il existe un anneau $A$ comme ci-dessus, de corps des fractions
$K$ tel que $\HH^N(\SP(K)_{\et},\ZZ/p)≠0$.
\emph{Supposons par l'absurde $N>d+r$}.
D'après les résultats du paragraphe \ref{Epp et maj grossière},
il existe alors une extension finie $L/K$ telle que 
$\HH_{\et}^N(L,\ZZ/p)=0$. Ce groupe est isomorphe à 
$\HH^N_{\et}(K,\mathrm{Ind}_K^L(\ZZ/p))$. 
La surjection canonique $\mathrm{Ind}_K^L(\ZZ/p)\surj \ZZ/p$ 
de modules galoisiens (donnée par la trace) induit ici,
puisqu'il n'y a pas de cohomologie en degré $≥N+1$,
une \emph{surjection} 
$$
0=\HH^N_{\et}(L,\ZZ/p)\surj \HH^{N}_{\et}(K,\ZZ/p)
$$
sur les groupes de cohomologie. Contradiction.

\section{Le cas d'un ouvert de $A[p^{-1}]$}\label{generalisation}

Dans cette section, on démontre le théorème suivant :

\begin{thm}\label{theoreme ouvert}
Soit $A$ un anneau hensélien excellent intègre de corps résiduel $k$ de
caractéristique $p>0$ et de corps des fractions de caractéristique nulle.
Alors, pour tout ouvert non vide affine $U\subset \SP(A[p^{-1}])$, on 
a :
$$
\mathrm{dim.coh.}_p(U_{\et})=\dim(A)+\dim_p(k).
$$
\end{thm}

Rappelons que pour tout topos $T$ on note $\mathrm{dim.coh.}_p(T)$
le plus grand entier $d$ tel que $\HH^{d}(T,\mc{F})≠0$ pour au moins un 
faisceau de $p$-torsion $\mc{F}$ (cf. \sga{4}{x}{§1} et 
\sga{4}{ix}{1.1}). 

On majore $\mathrm{dim.coh.}_p(U_{\et})$ en procédant
comme plus haut (cf. \ref{majoration ouvert}). 
Pour minorer $\mathrm{dim.coh.}_p(U_{\et})$, on se ramène au cas 
où $U$ est le point générique par une astuce de globalisation
d'une classe de cohomologie après un revêtement ramifié de
degré $2$ (cf. \ref{astuce Gabber}).
 
\subsection{Minoration}\label{minoration ouvert}

Soient $A,k$ comme en \ref{theoreme ouvert} ; posons $r=\dim(A)+\dim_p(k)$
que l'on suppose fini pour simplifier. On peut supposer $A$ normal car d'une
part $A$ est japonais et d'autre part, pour tout morphisme fini $A→A'$, on a
$\mathrm{dim.coh.}_p(A)≥\mathrm{dim.coh.}_p(A')$.
Soit $ξ$ un point maximal de $V(p)\subset \SP(A)$. Notons
$K^h_ξ$ le corps des fractions de l'hensélisé de l'anneau
de valuation discrète $A_ξ$. Par hypothèse,
$\dim_p(K^h_ξ)=r$ de sorte qu'il existe une extension finie
$K'/K^h_ξ$ telle que $μ_p\subset K'$ et $\HH^{r}(K',\ZZ/p)≠0$
(cf. \ref{coeff constants}).
L'extension (séparable) $K'/K^h_ξ$ peut être définie par un polynôme
irréductible unitaire à coefficients dans $K^h_ξ$. 
Par unicité de l'extension de la norme de $K^h_ξ$ à sa clôture 
séparable, le lemme de Krasner est applicable et l'on 
peut donc approximer les coefficients
de ce polynôme par des éléments de $K$ de sorte que l'extension
$K'/K^h_ξ$ soit définie sur $K$. Quitte à remplacer $A$ par sa normalisation
dans cette extension, on peut donc supposer que $\HH^{r}(K^h_ξ,μ_p^{⊗
r})≠0$. D'après les résultats de K. Katô (\cite{Galois@Kato}, théorème 1), 
ce groupe est engendré par les \emph{symboles} ; en particulier, il existe 
des éléments $φ₁,\dots,φ_r∈(K^h_ξ)^×$ tels que la classe
$$c=χ_{φ₁}\cup \cdots \cup χ_{φ_r}∈\HH^r(K^h_ξ,μ_p^{⊗r})$$
soit non nulle. Ci-dessus, $χ_{φ_i}$ est l'image de $φ_i$ par le 
composé $(K^h_ξ)^×→(K^h_ξ)^×/(K^h_ξ)^{×p}→\HH^{1}(K^h_ξ,μ_p)$.

Remarquons que l'anneau $A_ξ^h$ étant hensélien, il existe un entier $N$
tel que $$(\star)\ 1+\MM_ξ^N\subset (1+\MM_ξ)^p$$
dans cet anneau. Il en résulte immédiatement que l'on peut approximer
les $φ_i$ par des éléments de $K$ sans changer les $χ_{φ_i}$.

On va montrer dans la proposition suivante 
qu'il existe une classe de cohomologie
dans $\HH^r(\SP(A[p^{-1}]),ℱ)$, où $ℱ$ est un
$\ZZ/p$-faisceau, non nulle en restriction au corps des
fractions. La minoration désirée (pour tout ouvert affine
non vide $U$ comme dans l'énoncé du théorème) s'en déduit immédiatement.

\begin{prp2}[\cite{Algebrisation@Gabber}]\label{astuce Gabber}
Soient $A,\xi$ et $c$ comme ci-dessus. 
Pour tout voisinage ouvert $W\subset \SP(A)$ de $ξ$, il existe un morphisme
surjective, fini et plat de rang $2$, $π:U'→U:=\SP(A[p^{-1}])$,
tel que si l'on note $j'$ l'immersion ouverte $(W\cap U)':=π^{-1}(W\cap U)\hra U'$, 
il existe
également une classe $c'∈\HH^r(U',j'_!μ_p^{⊗r})$ telle que $π$ soit
(étale et) décomposé sur $U^h_ξ:=\SP(K^h_ξ)$ 
$$
U'×_U U^h_ξ\isononcan U^h_ξ \coprod U^h_ξ,
$$
de telle façon que $c'$ induise par restriction
$(c,-c)$ sur $U^h_ξ \coprod U^h_ξ$.
\end{prp2}

Il suffit de démontrer la proposition dans le cas particulier où $r=1$.
Remarquons tout d'abord à cette effet que si la conclusion de la proposition est vraie
pour un ouvert $W₀$, elle l'est pour tout ouvert $W$ le contenant.
Quitte à rétrécir $W$, on peut donc supposer 
les fonctions $φ_i$ inversibles sur $W\cap U$. 
Le produit (cf. p. ex. \sga{4,5}{Cycle}{1.2.4}) 
$$\HH^1((W \cap U)',μ_p)⊗\cdots ⊗\HH^1((W\cap U)',μ_p)⊗\HH^1(U',j'_!μ_p)→
\HH^r(U',j'_!μ_p^{⊗r})$$
$$χ_{φ_2}⊗\cdots⊗χ_{φ_r}⊗c'₁\mapsto c'$$
fournit une classe $c'$ à partir du cas $r=1$ et
des symboles associés aux $r-1$ dernières fonctions.

Supposons $W=\SP(A[h^{-1}])$. Quitte à diviser $φ$ par une puissance de $h^p$, 
on peut supposer que $φ^{-1}$ s'étend en une fonction
$ψ$ sur $U$, divisible par $h^2$. Enfin, quitte à multiplier $φ$ par une puissance 
convenable de $p^p$ (inversible sur $U$), on peut également 
supposer que $v_ξ(φ)≥N$. 

Considérons le revêtement $U'₀$ de $U$ défini par l'équation :
$$
f(X):=X^2-(ψ+2)X+1=0.
$$
Observons dès maintenant que la fonction $X$ est inversible sur $\OO_{U'₀}$.
Le polynôme \emph{unitaire} $g(X):=φ^2f(\frac{X}{φ})=X^2-(1+2φ)X+φ²$
(congru à $X(X-1)$ modulo $\MM_ξ$) 
satisfait aux conditions $g(1)∈(φ)$ (resp. $g(φ²)∈(φ^3)$)
et $g'(1)∈A_ξ^×$ (resp. $g'(φ²)∈A_ξ^×$) ; le polynôme
$f$ possède donc deux racines $x$ et $x'$ dans $K^h_ξ$
telles que $x ≡ φ \,\mathrm{mod}\,(φ^2)$ et $x' ≡ φ^{-1}\,\mathrm{mod}\,(1)$. 
Il existe donc $a∈A_ξ^h$ tel que $x=φ(1+aφ)$ ; par hypothèse sur $φ$ et $N$,
$(1+aφ)∈{K_ξ^{h×}}^{p}$, de sorte que $χ_x=χ_φ$ dans $\HH^1(K_ξ^h,μ_p)$.
Il en résulte immédiatement que $χ_{x'}=-χ_φ$. Soit $U'$ le normalisé
de $U'₀$ dans $U'[\frac{1}{h}]$ ; puisque par hypothèse $h^2$ divise
(dans $\OO_{U'₀}$) $ψ$, ainsi donc que son multiple $(X-1)²$, on a $h|X-1$ dans $\OO_{U'}$.
La fonction $X$ est donc une section sur $U'$ du faisceau $\GG_{m\,U',V(h)}:=\ker\big(\GG_{m\,U'}→
i_*\GG_{m\,V(h)}\big)$, où l'on note $i$ l'immersion fermée $V(h)\hra U'$.
Soit $j'$ l'immersion ouverte complémentaire. Puisque
$p$ est inversible sur $U'$, on a une suite exacte
$$
1→j'_!μ_p → \GG_{m\,U',V(h)}\sr{f\mapsto f^p}{→} \GG_{m\,U',V(h)}→1.
$$
qui étend la suite exacte usuelle de Kummer sur $K^h_ξ$.
L'image de $X$ par le morphisme cobord est la classe de cohomologie
$c'$ recherchée.

\subsection{Majoration}\label{majoration ouvert}

On procède par récurrence sur $\dim(A)$. Le cas de la dimension $1$ est dû 
à K. Katô. 
Le théorème de Lefschetz pour les morphismes affines de
type fini entre $\QQ$-schémas excellents (\sga{4}{xix}{6.1}),
nous ramène (grâce à l'hypothèse de récurrence)
au cas particulier où $U=\SP(A[p^{-1}])$.
La méthode de la trace (\sga{4}{ix}{5.6}) nous ramène à
montrer l'annulation de $\HH^n(\SP(A[p^{-1}]),\ZZ/p)$ pour
$n>\dim(A)+\dim_p(k)=:r$. On peut également supposer
$A$ normal et complet.
Soit $k₀$ le sous-corps parfait maximal de $k$. Comme
rappelé ci-dessus (§\ref{Epp et maj grossière}), il existe une extension
finie $W'/W(k₀)$ telle que le morphisme $\SP(A⊗_{W(k₀)} W')^{ν}→\SP(W')$
ait une fibre spéciale réduite. Procédant comme en §\ref{majoration cas corps},
le résultat pour $A[p^{-1}]⊗_{W(k₀)} W'$ se ramène par algébrisation au résultat
analogue dans le cas particulier où $A$ est 
l'hensélisation d'une algèbre de type fini sur
un anneau (local excellent) $B$ de dimension un de moins,
et que l'extension résiduelle pour $A/B$ est \emph{finie}. D'après
le théorème de Lefschetz affine (cf. \emph{op. cit.}), cela résulte
de l'hypothèse de récurrence.
La conclusion étant également valable pour les extensions finies de $W'$,
on a donc :
$$
\HH^n(A[p^{-1}]⊗_{K₀} \sur{K₀},\ZZ/p)=0\ ∀n>r,
$$
où l'on note $K₀=\Frac{\,W(k₀)}$.

Puisque $\mathrm{cd}_p(K₀)≤2$, il s'en suit que pour tout $\FF_p[G_{K₀}]$-module
$V$, et tout $n>r+2$, on a $\HH^n(A[p^{-1}],V)=0$.
Comme en §\ref{fin corps}, on observe que pour tout tel $V$, il 
existe une surjection $V'↠V$ où $V'$ est une somme directe d'induites
de représentations triviales de petits (\cad contenus dans $G_{\Frac{\,W'}}$)
sous-groupes ouverts. Pour un tel $V'$ on connaît le résultat
d'annulation pour $n>r$. L'annulation du $H^n$ pour $V$ se ramène
donc à l'annulation du $H^{n+1}$ pour $\ker(V'↠V)$. On obtient la majoration
désirée en procédant ainsi une fois de plus.

\section{Appendice : le théorème de structure de Cohen-Gabber}\label{appendice}

\begin{thm}[\cite{Conference-Illusie@Gabber}, lemme 8.1]\label{gabber}
Soit $A$ un anneau local complet noethérien \emph{réduit}, d'égale 
caractéristique $p>0$, équidimensionel
de dimension $d$ et de corps résiduel $k$.
Il existe un sous-anneau $A₀$ de $A$, isomorphe à $k[[t₁,\dots,t_d]]$,
tel que $A$ soit fini sur $A₀$, sans torsion et génériquement étale.
De plus, le morphisme $A₀→A$ induit un isomorphisme sur les corps résiduels.
\end{thm}

Ce résultat apparaît explicitement 
comme hypothèse (dans le cas intègre) dans \ega{$0_{\textsc{iv}}$}{21.9.5}.

La démonstration du théorème, tirée de \cite{Conference-Illusie@Gabber}, 
occupe le reste de cette section.

\subsection{}Soient $t_1,\dots,t_d∈A$ un système de paramètres,
$\{b_i\}_{i∈I}$ une $p$-base de $k=A/\MM_A$,
$\{β_i\}_{i∈I}$ des relèvements des $b_i$ dans $A$ (que
nous changerons par la suite), et $κ\subset A$ le corps de représentants 
correspondant (cf. \ac{ix}{2}{2}{}).

Notons $C$ l'ensemble des composantes irréductibles de $A$,
$\{℘_α\}_{α∈C}$ l'ensemble des idéaux premiers minimaux de $\SP(A)$, et 
$A_α:=A/℘_α$ l'anneau intègre de dimension $d$ correspondant
à la composante irréductible $α$. L'anneau $A$ étant réduit,
on a $(0)=\cap_{α} ℘_α$ ; c'est une décomposition primaire \emph{réduite} :
$\forall α, \cap_{β≠α} ℘_β \subsetneq ℘_α$.

Pour tout ensemble \emph{fini} $e\subset I$, posons $κ_e:=κ^p(β_i,\,i\notin
e)\subset κ$.
Les trois propriétés suivantes sont évidentes :
$$(\star)\ [κ:κ_e]<+∞,\, κ_{e\cup e'}\subset κ_{e}\cap κ_{e'}\ \textrm{et }
\cap_{e\subset I} κ_e = κ^p.$$

\subsection{}Pour simplifier les notations, fixons $α∈C$ et posons $B=A_α$,
$L$ son corps des fractions et $τ_i$ l'image de $t_i∈A$ dans $B$
par la surjection canonique.
Considérons les anneaux suivants :
$$R_κ=κ[[τ_1,\dots,τ_d]]\subset B,$$
$$L_κ=\Frac{\,R_κ}\subset L,$$
$$R_{κ,e}=κ_e[[τ_1^p,\dots,τ_d^p]]\subset R_κ,$$
$$L_{κ,e}=\Frac{\,R_{κ,e}}\subset L_κ ; $$

les morphismes d'inclusion sont finis (cf. p. ex.
\ega{$0_{\textrm{IV}}$}{19.8.8}, démonstration). Les
$τ_i$ ci-dessus sont analytiquement indépendants sur
$κ$ : $R_κ$ est un anneau de séries formelles.

D'après \cite{CRT@Matsumura}, §30, lemme 6, et
l'analogue de $(\star)$ pour les sous-corps $L_{κ,e}$ de 
$L_κ$, on a l'égalité 
$$
\rang_L Ω¹_{L/L_{κ,e}}=\rang_{L_κ} Ω¹_{L_κ/L_{κ,e}},
$$
dès que l'ensemble fini $e$ est suffisamment grand.

En particulier, on peut supposer cette égalité valable pour
chaque composante irréductible $α$.

Le terme de gauche est le rang (générique) du $B$-module
$Ω¹_{B/R_{κ,e}}$ ; remarquons que d'après \ega{$0_{\textrm{IV}}$}{21.9.4}, 
$Ω¹_{B/R_{κ,e}}$ s'identifie au module
$\chap{Ω}¹_{B/κ_e}$ des formes différentielles \emph{continues}.
Le terme de droite est le 
rang du $R_{κ,e}$-module libre $Ω¹_{R_{κ,e}/R_κ}$.
Ce dernier est égal à $d+\rang_κ Ω¹_{κ,κ_e}=d+|e|$
de sorte que l'on a :
\begin{equation}
\rang_B \chap{Ω}¹_{B/κ_e}=d+|e|
\end{equation}

\begin{lmm2}
Pour tout idéal non nul $I$ de $B$, l'ensemble
des $d(i)\otimes_B L$, pour $i∈I$, est une famille 
génératrice du $L$-espace vectoriel $\chap{Ω}¹_{B/κ_e}\otimes_B L$.
\end{lmm2}

\begin{proof}
En effet, si $i₀∈I$ est non nul et si l'on pose $v_0=d(i₀)$,
l'ensemble $d(Bi)=\{bv₀+i₀\mathrm{d}b,\,b∈B\}$ contient $v₀$ et est générateur
puisque l'ensemble des $\mathrm{d}b$ l'est.
\end{proof}

\subsection{}Supposons $e$ choisi comme ci-dessus
et $d>0$ sans quoi le théorème est trivial.
Identifions l'ensemble $C$ des composantes irréductibles de $\SP(A)$
à l'ensemble $\{1,\dots,c\}$, et notons pour tous $i∈e$ et $j∈\{1,\dots,c\}$,
$β_{i,j}$ l'image dans $A_j=A/℘_j$ de $β_i∈A$.
(Rappelons que les $β_i$ font partie d'une $p$-base de $κ\subset A$.)
Fixons $0≤j≤c-1$ et supposons qu'il existe des éléments $\{m_i\}_{i∈e}$ dans $\MM_A$
tels que les images des éléments $β_i+m_i$ dans chacun des
anneaux $A_1,\dots,A_j$ ait des différentielles
linéairement indépendantes (dans $Ω¹_{A₁/R_{κ,e}}\otimes_{A₁} \Frac{\,A₁}$,
\dots, $Ω¹_{A_j/R_{κ,e}}\otimes_{A_j} \Frac{\,A_j}$,). 
(Pour $j=0$, cette condition est vide.)
Vérifions qu'il en est de même pour $j+1$. On peut
supposer que les $m_i$ ($i∈e$) sont nuls ; nous le ferons pour 
simplifier les notations. Afin de ne pas altérer les
choix précédents sur les composantes $A_1,\dots,A_j$, on considère
l'idéal $℘_1\cap \cdots \cap ℘_j=\mathrm{Ker}(A→A_1×\cdots×A_j)$.
Comme rappelé plus haut, $℘_1\cap \cdots \cap ℘_j\subsetneq ℘_{j+1}$
de sorte que son image dans $B:=A/℘_{j+1}$ est un idéal non nul.
Si $j>0$, nous noterons $I$ son image ; si $j=0$, posons
$I=\MM_B$. D'après les résultats du paragraphe
précédent, $\rang_B \chap{Ω}¹_{B/κ_e}=d+|e|≥|e|$
et la famille $d(I)$ est génératrice dans $\chap{Ω}¹_{B/κ_e}\otimes_B L$
(où $L=\Frac{\,B}$). Il existe donc des éléments $m'_i∈I$, $i∈e$,
tels que les différentielles des éléments
$d(β_{i,j+1}+m_i')$, $i∈e$, soient linéairement indépendantes. Il 
suffit alors de relever les $m'_j$ dans $℘_1\cap \cdots \cap ℘_j\subset \MM_A$
si $j>0$, ou simplement dans $\MM_A$ si $j=0$,
pour obtenir les éléments souhaités. 

\subsection{}
D'après les résultats des deux paragraphes précédents, il existe un ensemble
fini $e$ tel que \emph{pour chaque composante irréductible} $B$ de $A$
on ait $\rang_B \chap{Ω}¹_{B/κ_e}=d+|e|$ et d'autre part
des éléments $β'_i$, $i∈e$, relevant les $b_i$,
dont les différentielles sont linéairement indépendantes
dans ces groupes. Si l'on considère le corps $κ':=κ^p(β_i, i\notin e; β'_i, i∈
e)=κ_e(β_i, i∈e)\subset A$, celui-ci s'envoie isomorphiquement sur $k$ par 
la surjection canonique $A→k=A/\MM_A$
et l'on a $$\rang_B \chap{Ω}¹_{B/κ'}=d$$ pour toute composante irréductible
$B$ de $A$. Pour simplifier les notations nous noterons ce nouveau
corps de représentants $κ$.

Le $A$-module $\chap{Ω}¹_{A/κ}$ étant de rang (générique) $d$ sur chaque
composante irréductible, on montre en procédant comme précédemment, 
qu'il existe des éléments $f_1,\dots,f_d$ de $A$ tels que les $d(f_i \ \mathrm{mod}
\, ℘_α)\otimes_{A_α} \Frac{\,A_α}$ forment une base de 
$\chap{Ω}¹_{A_α/κ} \otimes_{A_α} \Frac{\,A_α}$ pour chaque composante
irréductible $A_α$. Quitte à les multiplier par une puissance $p$-ième
d'un élément non nul appartenant à $\MM_A$, on voit que
l'on peut les supposer dans $\MM_A$. 
Rappelons que l'on a choisi un système de paramètres $t_1,\dots,t_d$
dans $A$, de sorte qu'en particulier, comme rappelé plus haut, 
le morphisme $A/k[[t_1,\dots,t_d]]$ est fini.

Posons, pour $i∈[1,d]$,
$$
t'_i:=t_i^p(1+f_i).
$$
Soit $A₀$ le sous-anneau $κ[[t'_1,\dots,t'_d]]$ de
$A$. Le morphisme $A/A₀$ est fini : cela résulte du 
fait que les éléments $1+f_i$ sont des unités de $A$.
Il est génériquement étale sur chaque composante
irréductible compte tenu de l'hypothèse
sur les éléments $f_i$ et de la formule
$$
d(t'_i)=t_i^p df_i.
$$

\bibliography{bib}
\bibliographystyle{smfalpha}

\end{document}